\documentclass[a4paper,11pt,onecolumn,twoside]{article}
\usepackage{mathrsfs}
\usepackage{mathrsfs}
\usepackage{amsfonts}
\usepackage{fancyhdr}
\usepackage{amsmath,amsfonts,amssymb,graphicx,amsthm}
\usepackage{cite}

\allowdisplaybreaks
\begin{document}

\title{\bf Nonlinear stability of a composite wave to the Cauchy problem of 1-D full compressible  Navier-Stokes-Allen-Cahn system}
\author{{\bf Dan Lei}\\
School of Mathematical Sciences, Anhui University, Hefei 230601, China\\[2mm]
{\bf Zhengzheng Chen}\thanks{Corresponding author.
{\it E-mail address: chenzz@ahu.edu.cn}}\\
Centre for Pure Mathematics, School of Mathematical Sciences, Anhui University,\\ Hefei 230601, China
}

\date{}

\vskip 0.2cm

\maketitle

\vskip 0.2cm \arraycolsep1.5pt
\newtheorem{Lemma}{Lemma}[section]
\newtheorem{Theorem}{Theorem}[section]
\newtheorem{Definition}{Definition}[section]
\newtheorem{Proposition}{Proposition}[section]
\newtheorem{Remark}{Remark}[section]
\newtheorem{Corollary}{Corollary}[section]

\begin{abstract}
 The compressible Navier-Stokes-Allen-Cahn system models the motion of a mixture of  two macroscopically immiscible viscous  compressible fluids. In this paper, we are concerned with  the large time behavior of solutions to the Cauchy problem of the one-dimensional  full compressible Navier-Stokes-Allen-Cahn system. If the Riemann problem of the  corresponding  Euler system admits a solution which is a linear combination  of 1-rarefaction wave and 3-rarefaction wave, we proved  that a global  strong solution to the  compressible Navier-Stokes-Allen-Cahn system exists uniquely and converges to the above composite wave as time goes to infinity, provided that  the adiabatic exponent $\gamma$ is closed to $1$.  Here  the initial perturbations except for the  temperature function of the fluid, and the strength of rarefaction waves can be arbitrarily large.  The proof is  given by an elementary energy method that takes  into account the effect of the phase field variable $\chi(t,x)$ and the complexity of nonlinear waves.
\bigbreak
\noindent

\bigbreak

\end{abstract}

\section{Introduction }
\setcounter{equation}{0}

This paper is concerned with a diffuse interface model proposed firstly by Blesgen \cite{Blesgen-1999}, namely, the compressible Navier-Stokes-Allen-Cahn (NSAC in short) system.  The model describes the motion of a mixture of  two macroscopically immiscible, viscous  compressible fluids which may undergo phase transitions.
In this model, the interface between the two-phase fluids is characterized   by a narrow transition layer. It is through this layer that the fluids may mix undergoing phase transitions. A phase field variable $\chi(t,x)\in[0,1]$ is introduced to indicate the location of the interface, and $\chi(t,x)$  is  governed by the modified Allen-Cahn equation. More precisely,  $\chi(t,x)$ takes the values $0$ and $1$ in the two absolute fluids respectively, and  $\chi(t,x)$ varies continuously in the interface of the two phases.  The interested readers are referred to   \cite{Abels-Feireisl,Feireisl,Rajagopal-2012,Kotschote-JEE,Kotschote-2017,Kotschote-2012-ARMA} for more physical background on such a model.

Chen et al. \cite{Shi-1} briefly described the compressible immiscible two-phase flow. This involves taking any one of the volume particles in the flow, and assuming that the mass of the components in the representative material volume $V$ is $m_i (i=1,2)$,  the mass concentration is $\chi_i=\frac{\rho_i}{\rho}$, the apparent mass density of the fluid $i$  is $\rho_i=\frac{m_i}{V}$. Then $\rho=\rho_1+\rho_2$ and $\chi=\chi_1-\chi_2$ represent the total density and the difference of the two components for the fluid mixture, respectively,  and the full  compressible NSAC system reads  as
\begin{eqnarray}\label{1.1}
\left\{\begin{array}{ll}
    \rho_t+\textnormal{div}(\rho {\bf u})=0,\\[2mm]
    (\rho{\bf u})_t+\textnormal{div}(\rho {\bf u}\otimes{\bf u})=\textnormal{div}\mathbb{T},\\[2mm]
    (\rho\chi)_t+\textnormal{div}(\rho\chi{\bf u})=-\mu,\\[2mm]
    \rho\mu=\rho\frac{\partial f}{\partial \chi}-\textnormal{div}\left(\rho\frac{\partial f}{\partial \nabla\chi}\right),\\[2mm]
    (\rho E)_t+\textnormal{div}(\rho E{\bf u})=\textnormal{div}\left(\mathbb{T}{\bf u}+\kappa\nabla\theta-\mu\frac{\partial f}{\partial \nabla\chi}\right),
\end{array}\right.
\end{eqnarray}
for $(t,\widetilde{{\bf x}})\in(0,\infty)\times\mathbb{R}^n$, where  $t$ and $\widetilde{{\bf x}}$  represent  the time variable and the spatial variable in the Eulerian coordinates, respectively, and $n$ is the spatial dimension.  The unknown functions $\rho(t,\widetilde{{\bf x}})$, ${\bf u}(t,\widetilde{{\bf x}})$, $\theta(t,\widetilde{{\bf x}})$ and $\chi(t,\widetilde{{\bf x}})$ denotes the total fluid density, the mean velocity, the absolute temperature of the fluid mixture and the concentration of one of the two fluids, respectively. $\mu(t,\widetilde{{\bf x}})$ is the chemical potential of the fluid, and $f(\rho,\chi,\nabla\chi)$ is the fluid-fluid interfacial free energy density.  $E$ denotes the total energy density of the mixture fluid, which is defined by
\begin{eqnarray}\label{1.2}
E:= e+\frac{1}{2}{\bf u}^2+ f(\rho,\chi,\nabla\chi),
\end{eqnarray}
here $e$ is the internal energy density, and  $f(\rho,\chi,\nabla\chi)$ is the interfacial free energy, which  is given by (see \cite{Rajagopal-2012})
\begin{eqnarray}\label{1.2-2}
f(\rho,\chi,\nabla\chi)=\frac{1}{4\epsilon}(1-\chi^2)^2+\frac{\epsilon}{2\rho}|\nabla\chi|^2
\end{eqnarray}
with the constant  $\epsilon>0$ being  the thickness of the diffuse interface.

The Cauchy stress-tensor $\mathbb{T}$ is represented by:
\begin{eqnarray}\label{1.3}
\mathbb{T}=2\widetilde{\nu}\mathbb{D}{\bf u}+\widetilde{\lambda}(\textnormal{div}{\bf u})\mathbb{I}-P(\rho,\theta,\nabla\chi)\mathbb{I}-\rho\nabla\chi\otimes\frac{\partial f}{\partial \nabla\chi},
\end{eqnarray}
where $P(\rho,\theta,\nabla\chi)$ is the pressure of fluid mixture, $\mathbb{I}$ is the unit matrix, $\mathbb{D}({\bf u})$ is the so-called deformation tensor satisfying
\begin{eqnarray}\label{1.5}
\mathbb{D}{\bf u}=\frac{1}{2}\left(\nabla{\bf u}+(\nabla{\bf u})^\top\right),
\end{eqnarray}
with the  superscript $\top$ denoting  the transpose, and all vectors here are column ones. Moreover,  $\widetilde{\nu},\widetilde{\lambda}>0$ are viscosity coefficients satisfying
\begin{eqnarray}\label{1.4}
\widetilde{\nu}>0,~\widetilde{\lambda}+\frac{2}{n}\widetilde{\nu}\geq 0,
\end{eqnarray}
and $\kappa>0$ denotes the heat conduction coefficient.

From the second law of thermodynamics, $P(\rho,\theta,\nabla\chi)$, $e(\rho,\theta)$ and $f(\rho,\chi,\nabla\chi)$ obey the following equation:
\begin{eqnarray}\label{1.6}
\theta ds=d(e+f)+Pd\left(\frac{1}{\rho}\right),
\end{eqnarray}
where $s$ is the entropy of the fluid mixture. From (\ref{1.6}), we know
\begin{eqnarray}\label{1.4-1}
\frac{\partial s}{\partial \theta}=\frac{1}{\theta}\frac{\partial (e+f)}{\partial \theta}, ~~~\frac{\partial s}{\partial \rho}=\frac{1}{\theta}\left(\frac{\partial (e+f)}{\partial \rho}-\frac{P}{\rho^2}\right).
\end{eqnarray}
Then, (\ref{1.4-1}) implies the following compatibility equation
\begin{eqnarray}\label{1.4-2}
P=\rho^2\frac{\partial (e+f)}{\partial \rho}+\theta\frac{\partial p}{\partial \theta}=\rho^2\frac{\partial e}{\partial \rho}-\frac{\epsilon}{2}|\nabla\chi|^2+\theta\frac{\partial p}{\partial \theta}.
\end{eqnarray}
Therefore, we consider the pressure $P(\rho,\theta,\nabla\chi)$ and the internal energy $e$ are given by
\begin{equation}\label{1.6-1}
P(\rho,\theta,\nabla\chi)=R\rho\theta-\frac{\epsilon}{2}|\nabla\chi|^2:=p(\rho,\theta)-\frac{\epsilon}{2}|\nabla\chi|^2,\quad e=C_v\theta+constant,
\end{equation}
where  $R>0$ is a constant, $C_v=\frac{R}{\gamma-1}$ is the specific heat  at constant volume, and  the constant $\gamma>1$ denotes the adiabatic exponent.

Throughout this paper, we assume that the viscosity coefficients $\widetilde{\nu},\widetilde{\lambda}>0$ and heat conduction coefficient $\kappa>0$ are all constants.
Then the full compressible NSAC system (\ref{1.1}) in $\mathbb{R}^n$ can be rewritten as
\begin{eqnarray}\label{1.7}
\left\{\begin{array}{ll}
    \rho_t+\textnormal{div}(\rho {\bf u})=0,\\[2mm]
    \rho{\bf u}_t+\rho ({\bf u}\cdot\nabla){\bf u}+\nabla p=\textnormal{div}\left(2\widetilde{\nu}\mathbb{D}{\bf u}+\widetilde{\lambda}(\textnormal{div}{\bf u})\mathbb{I}\right)+\textnormal{div}\left(\frac{\epsilon}{2}|\nabla\chi|^2-\epsilon\nabla\chi\otimes\nabla\chi\right),\\[2mm]
    \rho\chi_t+\rho{\bf u}\cdot\nabla\chi=-\mu,\\[2mm]
    \rho\mu=\frac{\rho}{\epsilon}(\chi^3-\chi)-\epsilon\Delta\chi,\\[2mm]
   C_v\rho\theta_t+\rho{\bf u}\cdot\nabla\theta+p\textnormal{div}{\bf u}=\textnormal{div}\left(\kappa\nabla\theta\right)+2\widetilde{\nu}|\mathbb{D}{\bf u}|^2+\widetilde{\lambda}(\textnormal{div}{\bf u})^2+\mu^2.
\end{array}\right.
\end{eqnarray}
In the one-dimensional whole space $\mathbb{R}$, the system (\ref{1.7})  is reduced to
\begin{eqnarray}\label{1.8}
\left\{\begin{array}{ll}
    \rho_t+(\rho u)_{\widetilde{x}}=0,\\[2mm]
    \rho u_t+\rho uu_{\widetilde{x}}+p_{\widetilde{x}}=\nu u_{\widetilde{x}\widetilde{x}}-\epsilon\left(\frac{\chi_{\widetilde{x}}^2}{2}\right)_{\widetilde{x}},\\[2mm]
    \rho\chi_t+ \rho u\chi_{\widetilde{x}}=-\mu,\\[2mm]
    \rho\mu=\frac{\rho}{\epsilon}(\chi^3-\chi)-\epsilon\chi_{{\widetilde{x}}{\widetilde{x}}},\\[2mm]
    C_v\rho\theta_t+C_v\rho u\theta_{\widetilde{x}}+pu_{\widetilde{x}}=\kappa\theta_{\widetilde{x}\widetilde{x}}
    +\nu u_{\widetilde{x}}^2+\mu^2
\end{array}\right.
\end{eqnarray}
for $(t,\widetilde{x})\in\mathbb{R^+}\times\mathbb{R}$, where $\nu=2\widetilde{\nu}+\widetilde{\lambda}>0$.\\

Let us recall some known results about the compressible NSAC system. Many researchers have dedicated their efforts to the compressible NSAC system in a bounded domain. More specifically, Feireisl et al. \cite{Feireisl}, and Chen et al. \cite{Chen-Wen-Zhu,Chen-Ji-Wen-Zhu} investigated the existence of global-in-time weak solutions to the 3D isentropic system. Axmann and Mucha\cite{Axmann-Mucha-2016} obtained the existence of weak solutions to the 1D stationary isentropic system with slip boundary conditions. Li et al. \cite{Li-Yan-Ding-Chen} established the global existence of weak solutions to the 1D isentropic system with free boundary, provided that the density connects with vacuum continuously. The local existence and uniqueness of strong solutions to the full system in $\mathbb{R}^n$ have been investigated by Kotschote \cite{Kotschote-2012-ARMA}. Chen et al. \cite{Chen-He-Huang-Shi} established the existence and uniqueness of local strong solutions to the isentropic compressible model with  Navier and relaxation boundary conditions in a bounded domain of $\mathbb{R}^3$. Ding et al. \cite{Ding-2013}, Chen and Guo \cite{Chen-guo} proved the global existence of solutions to the Cauchy and initial boundary value problem of the 1D isentropic system for initial data without and with vacuum states. Chen and Zhu \cite{Zhu-ZAMP} established the blow-up criteria and global existence of strong solutions to the 1D isentropic system. Ding et al. \cite{Ding-2019,Ding-2024} also obtained the global existence of strong solutions to the free boundary value problems of the 1D isentropic system. Shi et al. \cite{Shi-1}, and Yan, Ding and Li \cite{Yan-Ding-Li-2022} closely investigated the global well-posedness of strong solutions to the 1D full system. Dai et al. \cite{Dai-Ding-Li} obtained the global existence of strong solutions to the 1D full system. Moreover, Song et al. \cite{Song-Zhang-Wang} proved the existence and uniqueness of the time-periodic solution to the compressible isentropic system with a time periodic external force in a periodic domain of $\mathbb{R}^n$. Note that the initial data in the above literatures \cite{Feireisl,Chen-Wen-Zhu,Chen-Ji-Wen-Zhu,Chen-He-Huang-Shi,Li-Yan-Ding-Chen,Dai-Ding-Li,Kotschote-2012-ARMA,Axmann-Mucha-2016,Ding-2013,Chen-guo,Zhu-ZAMP,Ding-2019,Ding-2024,Shi-1,Yan-Ding-Li-2022} can be arbitrarily large.

The compressible NSAC system in some unbounded domain has also attracted a considerable number of researchers. As far as we know, there have been a large number of papers on the large-time behavior of solutions to the initial/initial-boundary value problem of the compressible NSAC equations toward the viscous version of the basic wave patterns, namely, shock wave, rarefaction wave, contact discontinuity, stationary wave and their compositions. More precisely, Luo et al. \cite{Luo-Yin-Zhu,Luo-43,Luo-zhu-3MAS}, and Chen, Huang and Shi \cite{Chen-Huang-Shi,Chen-Huang-Shi-2024} proved the nonlinear stability of some basic waves, such as the viscous shock wave, the  rarefaction wave, the contact discontinuity and their linear combinations to the Cauchy problem for the 1D isentropic/non-isentropic system under small initial perturbations. Luo \cite{Luo-2022}, Yin and Zhu \cite{Yin-Zhu-JDE-2019} investigated the asymptotic stability of stationary solution and superposition of a stationary solution and a rarefaction wave for inflow problem of the 1D non-isentropic system under small initial perturbations. Lei and Chen \cite{Lei-Chen-} proved the asymptotic stability of rarefaction wave to the impermeable wall problem for the 1D isentropic system with density and phase field variable dependent viscosity, provided that the initial perturbations are partially large. Afterward, Chen, Lei and Yin \cite{Stationary-solutions} obtained the existence, stability and convergence rate of stationary solutions to the outflow problem for the 1D full system under small initial perturbations. Moreover, the global existence and optimal time decay rates of classical solutions to the Cauchy problem for the 3D isentropic system \cite{Zhao,Chen-Hong-Shi} and non-isentropic system \cite{Chen-He-Tang} were also established under small initial perturbations. Chen, Tang and Zhang \cite{Chen-Tang-Zhang} studied the global existence and the pointwise behavior of strong solutions to the Cauchy problem of the 3D compressible isentropic system under the assumption that the initial perturbation is small and decays spatially.

However, to our knowledge, there are few results about the large-time behavior of solutions to the initial value problem for the 1D compressible isentropic/non-isentropic NSAC system in some unbounded domain with arbitrarily large initial perturbation. Chen et al. \cite{Chen-DCDS-B,Lei-Chen-JMAA} studied the global existence and large-time behavior of  strong solutions to the Cauchy problem of the 1D isentropic model with density/phase filed variable dependent viscosity and partially large initial data. Chen and Lei \cite{Chen-Lei} proved the global existence of strong solutions to the initial/initial-boundary value problem for the 1D isentropic system with arbitrarily large initial data. Thus a natural central question is: {\it Is it possible to achieve the global nonlinear stability of the composite wave to the Cauchy problem for the 1D full compressible NSAC system in an unbounded domain with arbitrarily large initial data?} To address this question, we focus on a specific scenario in the present work: this is the main intent of this paper, and we shall consider the existence of global-in-time rarefaction wave solutions to the 1D non-isentropic NSAC system (\ref{1.8}) with arbitrarily large initial data in unbounded domain. Here, we will give a partial  positive answer to this important problem.

In this paper, we are concerned with the global nonlinear stability of the composite wave consisting of two rarefaction waves for the Cauchy problem to the 1D full compressible NSAC system (\ref{1.8}) with arbitrarily large initial data. Now we turn to our original problem (\ref{1.8}). To study such a problem more easily, let's transform it into the Lagrangian coordinates $(t,x)$. Namely, set
\begin{eqnarray*}
x=\int_0^{\widetilde{x}}\rho(\tau,z)dz-\int_0^\tau(\rho u)(\xi,0)d\xi,~~\qquad t=\tau,
\end{eqnarray*}
then the 1D non-isentropic compressible NSAC system (\ref{1.8}) in the Lagrangian coordinates can be rewritten as
\begin{eqnarray}\label{1.9}
\left\{\begin{array}{ll}
    v_t-u_x=0,\\[2mm]
    u_t+p_x=\nu\displaystyle\left(\frac{u_x}{v}\right)_x-\frac{\epsilon}{2}\left(\frac{\chi_x^2}{v^2}\right)_x,\\[2mm]
    \chi_t=-v\mu,\\[2mm]
    \mu=-\epsilon\displaystyle\left(\frac{\chi_x}{v}\right)_x+\frac{1}{\epsilon}(\chi^3-\chi),\\[2mm]
    C_v\theta_t+pu_x=\kappa\displaystyle\left(\frac{\theta_x}{v}\right)_x+\nu\frac{u_x^2}{v}+v\mu^2,
\end{array}\right.\,\,\,t>0,~ x\in\mathbb{R}
\end{eqnarray}
with the prescribed initial data and far field conditions:
\begin{eqnarray}\label{1.9-1}
&&(v,u,\theta,\chi)(t,x)\big|_{t=0}=(v_{0},u_{0},\theta_0,\chi_0)(x),\,\, x\in\mathbb{R},~\inf_{x\in\mathbb{R}}v_0(x)>0,~\inf_{x\in\mathbb{R}}\theta_0(x)>0,\nonumber\\
&&\lim_{x\rightarrow\pm\infty}(v,u,\theta,\chi)(t,x)=(v_{\pm},u_{\pm},\theta_{\pm},\chi_{\pm}),~~t>0,
\end{eqnarray}
where $t$ and $x$ denote the time and the spatial variable in the Lagrangian coordinates, respectively. $v_{\pm}>0,~u_{\pm},~\theta_{\pm}>0$ and $\chi_{\pm}$ are constants. We assume, as usual in thermodynamics, that given any two of the five thermodynamical variables, $v,~p,~e,~\theta$~and $s$, the remaining three variables can expressed by these two given variables. Here and after, we only assume that  the thickness of the diffuse interface is $\epsilon=1$ for the present study.

Notice that, the second law of thermodynamics asserts that
\begin{eqnarray}\label{1.10}
\theta ds=d(e+f)+Pdv,
\end{eqnarray}
where $P=p-\frac{1}{2}\chi_x^2$, from which, if we choose $(v,\theta)$, or $(v,s)$  as independent variables and write $(p,e,s)=(p,e,s)(v,\theta)$, or $(p,e,\theta)=(\widetilde{p},\widetilde{e},\widetilde{\theta})(v,s)$, respectively, then we can deduce that:
\begin{eqnarray}\label{1.11}
\left\{\begin{array}{ll}
    s_v(v,\theta)=p_{\theta}(v,\theta),\\[2mm]
    s_\theta(v,\theta)=\frac{e_{\theta}(v,\theta)}{\theta},\\[2mm]
    e_v(v,\theta)=\theta p_{\theta}(v,\theta)-p(v,\theta),
\end{array}\right.
\end{eqnarray}

\begin{eqnarray}\label{1.12}
\left\{\begin{array}{ll}
    \widetilde{e}_v(v,s)=-p(v,\theta),\,\,\,~~~~~~~~~~~~~~~\widetilde{e}_s(v,s)=\theta\\[2mm]
    \widetilde{p}_v(v,s)=p_v(v,\theta)-\frac{\theta (p_{\theta}(v,\theta))^2}{e_\theta(v,\theta)},\,\,\,\widetilde{p}_s(v,s)=\frac{\theta p_{\theta}(v,\theta)}{e_\theta(v,\theta)},\\[2mm]
    \widetilde{\theta}_v(v,s)=-\frac{\theta p_{\theta}(v,\theta)}{e_\theta(v,\theta)},\,\,\,~~~~~~~~~~~~~~\widetilde{\theta}_s(v,s)=\frac{\theta}{e_{\theta}(v,\theta)},
\end{array}\right.
\end{eqnarray}

In this paper, we aim to establish the nonlinear stability of expansion waves for the system (\ref{1.9}). To this end, it is more convenient to work with the equations for the entropy $s$ and the absolute temperature $\theta$, that is
\begin{eqnarray}\label{1.13}
s_t=\kappa\left(\frac{\theta_x}{v\theta}\right)_x+\kappa\frac{\theta_x^2}{v\theta^2}+\nu\frac{u_x^2}{v\theta}+\frac{v\mu^2}{\theta},
\end{eqnarray}
where we used $(\ref{1.9})_5$, $(\ref{1.12})_5$ and $(\ref{1.12})_6$.

Indeed, for smooth solutions, system $(\ref{1.9})$ is equivalent to system $(\ref{1.9})_1$-$(\ref{1.9})_4$, (\ref{1.13}). In what follows, we will consider the Cauchy problem of system $(\ref{1.9})_1$-$(\ref{1.9})_4$, (\ref{1.13}).

At this point, we have to work out the relationships among $\theta,v, s$.

From $(\ref{1.6-1})_1$, $(\ref{1.6-1})_2$ and $(\ref{1.12})_6$, it is easy to that
\begin{eqnarray}\label{1.13-1}
p(v,\theta)=R\frac{\theta}{v},~~~~\theta=C(v)\textnormal{exp} \left(\frac{\gamma-1}{R}s\right).
\end{eqnarray}
Next, we need to determine the expression for $C(v)$.
Due to $(\ref{1.12})_1$, $(\ref{1.12})_2$ and  $(\ref{1.13-1})_1$, we see that $\widetilde{e}_{sv}=\widetilde{e}_{vs}$, that is, $\theta_v=-p_s=-\frac{R}{v}\theta_s$. Then combining the identity with $(\ref{1.13-1})_2$, we can get
\begin{eqnarray}\label{1.13-2}
C(v)=v^{-\gamma+1},
\end{eqnarray}
which implies that
\begin{eqnarray}\label{1.13-3}
\theta=v^{-\gamma+1}\textnormal{exp} \left(\frac{\gamma-1}{R}s\right).
\end{eqnarray}
Inserting (\ref{1.13-3}) into $(\ref{1.13-1})_1$, the pressure $p(v,\theta)$  satisfies the following equation of the state:
\begin{eqnarray}\label{1.15}
p(v,\theta)=R\frac{\theta}{v}=Av^{-\gamma} \textnormal{exp} \left(\frac{\gamma-1}{R}s\right),
\end{eqnarray}
where the specific gas constants $A$ and $R$ are positive constants, and $\gamma>1$ is the adiabatic constant.

The initial data for system $(\ref{1.9})_1$-$(\ref{1.9})_4$, (\ref{1.13}) are given by:
\begin{eqnarray}\label{1.16}
&&(v,u,s,\chi)(t,x)\big|_{t=0}=(v_{0},u_{0},s_0,\chi_0)(x),\,\, x\in\mathbb{R},~\inf_{x\in\mathbb{R}}v_0(x)>0,
\end{eqnarray}
and the far fields $x=\pm\infty$ data is given by
\begin{eqnarray}\label{1.17}
&&\lim_{x\rightarrow\pm\infty}(v,u,s,\chi)(t,x)=(v_{\pm},u_{\pm},s_{\pm},\chi_{\pm}),~~t>0,
\end{eqnarray}
where $v_{\pm}>0,~u_{\pm},~s_{\pm}$ and $\chi_{\pm}$ are constants, and the following compatibility conditions hold:
\begin{eqnarray}\label{1.18}
\rho_0(\pm\infty)=\rho_{\pm},~u_0(\pm\infty)=u_{\pm},~s_0(\pm\infty)=s_\pm,~\chi_0(\pm\infty)=\chi_{\pm}.
\end{eqnarray}
Here $s_0(x)=\frac{R}{\gamma-1}\ln\left(\frac{R}{A}\theta_0(x)v_0^{\gamma-1}(x)\right)$ and $s_\pm=\frac{R}{\gamma-1}\ln\left(\frac{R}{A}\theta_\pm v_\pm^{\gamma-1}\right)$. Since we consider the expansion waves to $(\ref{1.9})_1$-$(\ref{1.9})_4$, (\ref{1.13}), we assume that $s_+=s_-=\overline{s}$ and $\chi_+=\chi_-=1$ in the rest of this paper.

To study the large time behavior of solutions to the Cauchy problem $(\ref{1.9})_1$-$(\ref{1.9})_4$, (\ref{1.13}), it is noted that in the simplified case of the phase field variable $\chi\equiv1$, the large time behavior of this Cauchy problem corresponds to that of   the  classical compressible Navier-Stokes system given blow:
\begin{eqnarray}\label{1.19}
\left\{\begin{array}{ll}
    v_t=u_x,~~\,\,\,t>0,~ x\in\mathbb{R},\\[2mm]
    u_t+p(v,\theta)_x=\nu\left(\frac{u_x}{v}\right)_x,~~\,\,\,t>0,~ x\in\mathbb{R},\\[2mm]
    C_v\theta_t+p(v,\theta)u_x=\kappa\left(\frac{\theta_x}{v}\right)_x+\frac{\nu u_x^2}{v},~~\,\,\,t>0,~ x\in\mathbb{R},\\[2mm]
    (v,u,\theta)(t,x=\pm\infty)=(v_{\pm},u_{\pm},\theta_{\pm}),~~t>0,\\[2mm]
    (v,u,\theta)(0,x)=(v_{0},u_{0},\theta_0)(x)\rightarrow(v_{\pm},u_{\pm},\theta_{\pm}),~~\textnormal{as}~x\rightarrow\pm\infty.
\end{array}\right.
\end{eqnarray}

\subsection{Euler equations and composite wave}
For the expansion wave, the right-hand side of (\ref{1.19}) decays faster than each individual term on the left-hand side. Therefore, the compressible Navier-Stokes system (\ref{1.19}) may be approximated time-asymptotically by the Riemann problem of the compressible Euler equations:
\begin{eqnarray}\label{1.20}
\left\{\begin{array}{ll}
    v_t=u_x,\\[2mm]
    u_t+\widetilde{p}(v,s)_x=0,\\[2mm]
    s_t=0,
\end{array}\right.~~\,\,\,t>0,~ x\in\mathbb{R}
\end{eqnarray}
with Riemann initial data:
\begin{eqnarray}\label{1.21}
(v,u,s)(t,x)\big|_{t=0}=\left(v_0^R,u_0^R,s_0^R\right)(x)=\left\{\begin{array}{ll}
    (v_{-},u_{-},s_{-}),~~x<0,\\[2mm]
    (v_{+},u_{+},s_{+}),~~x>0.
\end{array}\right.
\end{eqnarray}

It is well-known that the Riemann problem (\ref{1.20})-(\ref{1.21}) admits a unique global weak entropy solution $\left(V^R,U^R,S^R\right)\left(\frac{x}{t}\right)$, which consists of a rarefaction wave of the first family, denoted by $\left(V_1^R\left(\frac{x}{t}\right),U_1^R\left(\frac{x}{t}\right),\overline{s}\right)$, and another of the third family, denoted by $\bigl(V_3^R\left(\frac{x}{t}\right),U_3^R\left(\frac{x}{t}\right),
\overline{s}\bigr)$. That is, there exists a unique constant state $(v_m,u_m)\in\mathbb{R}^2(v_m>0)$ such that $(v_m,u_m)\in\mathbb{R}_1(v_-,u_-)$ and $(v_+,u_+)\in\mathbb{R}_3(v_m,u_m)$. Here
\begin{equation*}
\mathbb{R}_1(v_-,u_-,\overline{s})=\left\{(v,u,s)\big| u=u_-+\int_{v_-}^v\sqrt{A\gamma z^{-\gamma-1}\exp((\gamma-1)\overline{s})}dz,u\geq u_-,s=\overline{s}\right\},
\end{equation*}
and
\begin{equation*}
\mathbb{R}_3(v_m,u_m,\overline{s})=\left\{(v,u,s)\big| u=u_m-\int_{v_m}^v\sqrt{A\gamma z^{-\gamma-1}\exp((\gamma-1)\overline{s})}dz,u\geq u_m,s=\overline{s}\right\}.
\end{equation*}

In other words, the unique weak entropy solution $\left(V^R,U^R,S^R\right)\left(\frac{x}{t}\right)$ to the Riemann problem (\ref{1.20})-(\ref{1.21}) is given by
\begin{eqnarray}\label{1.22}
\left(V^R,U^R,S^R\right)\left(\frac{x}{t}\right)=\left(V_1^R\left(\frac{x}{t}\right)+V_3^R\left(\frac{x}{t}\right)-v_m,U_1^R\left(\frac{x}{t}\right)+U_3^R\left(\frac{x}{t}\right)-u_m,\overline{s}\right),
\end{eqnarray}
where $\left(V_i^R,U_i^R,S^R\right)\left(\frac{x}{t}\right)(i=1,3)$ are determined by the following equations:
\begin{equation}\label{1.23}
\left(V_1^R,U_1^R\right)\left(\frac{x}{t}\right)=\left\{\begin{array}{ll}
    (v_-,u_-),~~~~~~~~~~~~~~~~~~~~~~~~~~~~~~~~~~~~-\infty\leq\frac{x}{t}\leq \lambda_1(v_-,\overline{s}),\\[2mm]
    \left(\lambda_1^{-1}\left(\frac{x}{t}\right),u_--\int_{v_-}^{V_1^R\left(\frac{x}{t}\right)}\lambda_1(z,\overline{s})dz\right),~~\lambda_1(v_-,\overline{s})\leq\frac{x}{t}\leq \lambda_1(v_m,\overline{s}),\\[2mm]
    (v_m,u_m),~~~~~~~~~~~~~~~~~~~~~~~~~~~~~~~~~~~~\lambda_1(v_m,\overline{s})\leq\frac{x}{t}\leq+\infty,
\end{array}\right.
\end{equation}
and
\begin{equation}\label{1.24}
\left(V_3^R,U_3^R\right)\left(\frac{x}{t}\right)=\left\{\begin{array}{ll}
    (v_m,u_m),~~~~~~~~~~~~~~~~~~~~~~~~~~~~~~~~~~~~-\infty\leq\frac{x}{t}\leq \lambda_3(v_m,\overline{s}),\\[2mm]
    \left(\lambda_3^{-1}\left(\frac{x}{t}\right),u_m-\int_{v_m}^{V_3^R\left(\frac{x}{t}\right)}\lambda_3(z,\overline{s})dz\right),~~\lambda_3(v_m,\overline{s})\leq\frac{x}{t}\leq \lambda_3(v_+,\overline{s}),\\[2mm]
    (v_+,u_+),~~~~~~~~~~~~~~~~~~~~~~~~~~~~~~~~~~~~\lambda_3(v_+,\overline{s})\leq\frac{x}{t}\leq+\infty
\end{array}\right.
\end{equation}
with $\lambda_1=\lambda_1(v,\overline{s})=-\sqrt{-\widetilde{p}_v(v,\overline{s})}$ and $\lambda_3=\lambda_3(v,\overline{s})=\sqrt{-\widetilde{p}_v(v,\overline{s})}$.

To analyze the nonlinear stability of the expansion waves, we first construct a smooth approximation to the above Riemann solution (\ref{1.22}),  following the approach in \cite{Matsumura-1986}. For this purpose, we begain with the Riemann problem for the typical Burgers equation:
\begin{eqnarray}\label{1.25}
\left\{\begin{array}{ll}
    \omega_{it}^R+\omega_i^R\omega_{ix}^R=0,~~t>0,~x\in\mathbb{R},\\[2mm]
    \omega_i^R(t,x)\big|_{t=0}=\omega_{i0}^R(x)=\left\{\begin{array}{ll}
    \omega_{i-},~~x<0,\\[2mm]
    \omega_{i+},~~x>0,
    \end{array}\right.~i=1,3,
\end{array}\right.
\end{eqnarray}
where
\begin{eqnarray}\label{1.26}
\left\{\begin{array}{ll}
      \omega_{1-}=\lambda_1(v_-,\overline{s})=-\sqrt{-\widetilde{p}_v(v_-,\overline{s})},~~\omega_{1+}=\lambda_1(v_m,\overline{s})=-\sqrt{-\widetilde{p}_v(v_m,\overline{s})},\\[2mm]
      \omega_{3-}=\lambda_3(v_m,\overline{s})=\sqrt{-\widetilde{p}_v(v_m,\overline{s})},~~~~\omega_{3+}=\lambda_3(v_+,\overline{s})=\sqrt{-\widetilde{p}_v(v_+,\overline{s})}.
\end{array}\right.
\end{eqnarray}
Then the approximate rarefaction wave of $\left(V^R,U^R,
S^R\right)\left(\frac{x}{t}\right)$ can be defined by
\begin{eqnarray}\label{1.29}
&&\left(V,U,S\right)(t,x)
=\left(V_1(t+t_0,x)+V_3(t+t_0,x)-v_m,U_1(t+t_0,x)+U_3(t+t_0,x)-u_m,\overline{s}\right),~~
\end{eqnarray}
where $t_0=\frac{1}{\varepsilon^2}$, and $(V_i,U_i)(t,x)(i=1,3)$ satisfy
\begin{eqnarray}\label{1.30}
&&\lambda_i(V_i(t,x),\overline{s})=\omega_i(t,x),~~i=1,3,\nonumber\\
&&\lambda_1(v,s)=-\sqrt{-\widetilde{p}_v(v,s)},~~~\lambda_3(v,s)=\sqrt{-\widetilde{p}_v(v,s)},\\
&&U_1(t,x)=u_-+\int_{v_-}^{V_1(t,x)}\sqrt{-\widetilde{p}_z(z,s)}dz,~~~U_3(t,x)=u_m-\int_{v_m}^{V_3(t,x)}\sqrt{-\widetilde{p}_z(z,s)}dz.\nonumber
\end{eqnarray}
Correspondingly, we can define $\Theta_i(t,x)(i=1,3)$ as
\begin{eqnarray}\label{1.31}
\Theta_i(t,x)=\widetilde{\theta}(V_i(t,x),\overline{s}),~i=1,3,~\textnormal{with}~\theta_m=\widetilde{\theta}(v_m,\overline{s}).
\end{eqnarray}

\subsection{Main result}

With the above preparations,  our main result can be stated as follows.
\begin{Theorem} For any given left end state $(v_-,u_-,\theta_-)$, suppose that the right end state $(v_+,u_+,\theta_+)$ is connected to the left end state $(v_-,u_-,\theta_-)$ by $(V^R,U^R,\Theta^R,S^R)(t,x)$, which is the solution to the Riemann problem of the compressible Euler system (\ref{1.20})-(\ref{1.21}) given by (\ref{1.22}). Moreover, $(V,U,\Theta,S)(t,x)$ is a smooth approximation of $(V^R,U^R,\Theta^R,
S^R)(t,x)$ constructed by (\ref{1.29}). Assume that
\begin{itemize}
\item[$\bullet$] $N_{0i}:=\left\|\left(v_0(x)-V(0,x),u_0(x)-U(0,x),\frac{\theta_0(x)-\Theta(0,x)}{\sqrt{\gamma-1}},\chi_0(x)-1\right)\right\|_i$ (i=,1,2) is bounded by some positive constant independent of $\gamma-1$;\\
\item[$\bullet$] There exist $(\gamma-1)$-independent positive constants $0<\underline{V}<1,~\overline{V}>1$, $0<\underline{\Theta}<1,~\overline{\Theta}>1$ such that
    \begin{equation*}
    2\underline{V}\leq v_0(x), V(t,x)\leq \frac{1}{2}\overline{V}, ~~2\underline{\Theta}\leq \theta_0(x), \Theta(t,x)\leq \frac{1}{2}\overline{\Theta}, ~~\forall (t,x)\in\mathbb{R}^+\times\mathbb{R};
    \end{equation*}
\item[$\bullet$] The constants $v_\pm,u_\pm$, and $\theta_\pm$ do not depend on $\gamma-1$;
\item[$\bullet$] There exists a sufficiently small positive constant $\varepsilon_0$ which depends only on $N_{0i}$, $\underline{V},\overline{V},\underline{\Theta}$ and $\overline{\Theta}$ such that
    \begin{equation*}
    0<\gamma-1\leq\varepsilon_0.
    \end{equation*}
\end{itemize}

Then the Cauchy problem $(\ref{1.9})_1$-$(\ref{1.9})_4$, (\ref{1.13}) admits a unique global smooth solution $(v,u,\theta,s,\chi)\\(t,x)$ which satisfies
\begin{equation}\label{1.34}
C_0^{-1}\leq v(t,x)\leq C_0,~~\underline{\Theta}\leq\theta(t,x)\leq\overline{\Theta},~0\leq\chi(t,x)\leq 1,~\forall (t,x)\in\mathbb{R}^+\times\mathbb{R},
\end{equation}
and
\begin{eqnarray}\label{1.35}
&&\sup_{0\leq t\leq+\infty}\left\|\left(v-V,u-U,\frac{\theta-\Theta}{\sqrt{\gamma-1}},s-\overline{s},\chi-1\right)(t)\right\|_2^2\nonumber\\
&&+\int_0^{+\infty}\left(\|(v-V)_x(\tau)\|_1^2+\left\|(u-U,\theta-\Theta)_x(\tau)\right\|_2^2+\|(\chi-1)(\tau)\|_3^2\right)d\tau\\
\leq&&C_1\left\|\left(v_0(\cdot)-V(0,\cdot),u_0(\cdot)-U(0,\cdot),\frac{\theta_0(\cdot)-\Theta(0,\cdot)}{\sqrt{\gamma-1}},s_0(\cdot)-\overline{s},\chi_0(\cdot)-1\right)\right\|_2^2.\nonumber
\end{eqnarray}
Here $C_0$ is a positive constant depending only on $\underline{V},~\overline{V},~\underline{\Theta},~\overline{\Theta}$, $\delta$,~$\|v_0(\cdot)-V(0,\cdot)\|,~\|u_0(\cdot)-U(0,\cdot)\|,\\
\left\|\frac{\theta_0(\cdot)-\Theta(0,\cdot)}{\sqrt{\gamma-1}}\right\|$~and~$\|\chi_0(\cdot)-1\|_1$, and $C_1$ is a positive constant depending only on $\underline{V},~\overline{V},~\underline{\Theta},~\overline{\Theta}$, $\delta$,~$\|v_0(\cdot)-V(0,\cdot)\|_2,~\|u_0(\cdot)-U(0,\cdot)\|_2,~
\left\|\frac{\theta_0(\cdot)-\Theta(0,\cdot)}{\sqrt{\gamma-1}}\right\|_2$ and $\|\chi_0(\cdot)-1\|_2$.

Furthermore,
\begin{equation}\label{1.36}
\lim_{x\rightarrow+\infty}\sup_{x\in\mathbb{R}}\left|(v-V,u-U,\theta-\Theta,s-\overline{s},\chi-1)(t,x)\right|=0.
\end{equation}
\end{Theorem}
 Some remarks on Theorem 1.1 are listed as follows.
\begin{Remark}
In Theorem 1.1, $\gamma-1$ needs to be small, but the initial perturbations for $\|v_0-V(0,\cdot)\|_2,~\|u_0-U(0,\cdot)\|_2$, $\left\|\frac{\theta_0-\Theta(0,\cdot)}{\sqrt{\gamma-1}}\right\|_2$ and $\|\chi_{0}-1\|_2$, and the strength of the rarefaction waves $\delta=|v_+-v_-|+|u_+-u_-|$ can be arbitrarily large. Thus in this sense, we obtain the time-asymptotic stability of composite wave to the Cauchy problem for the 1-D compressible non-isentropic NSAC system with arbitrarily large initial perturbation on the whole space $\mathbb{R}$. Note that for the previous results \cite{Chen-DCDS-B,Lei-Chen-JMAA,Lei-Chen-,Chen-Lei} on the 1-D compressible isentropic NSAC system in some unbounded domain with arbitrarily/partially large initial data, our results are more innovative and challenging than ever.
\end{Remark}
\begin{Remark}
From the proof of Theorem 1.1, $\gamma-1$ is only assumed to be sufficiently small such that
\begin{eqnarray*}
(\gamma-1)\left(\left\|\left(v_0-V(0,\cdot),u_0-U(0,\cdot),\frac{\theta_0-\Theta(0,\cdot)}{\sqrt{\gamma-1}}\right)\right\|^2_1+\|\chi_{0}-1\|_2^2+1\right)<C
\end{eqnarray*}
holds for some positive constant $C$ independent of $\gamma-1$ [see (\ref{2.65})]. Noted that
\begin{eqnarray*}
\left\|\frac{\theta_0-\Theta(0,\cdot)}{\sqrt{\gamma-1}}\right\|_1\leq C(\underline{V},\overline{V})\sqrt{\gamma-1}\left\|e^{(\gamma-1)s_0(x)}\right\|_{L^{\infty}(\mathbb{R})}\|(v_0-V(0,\cdot),s_0(x))\|_1,
\end{eqnarray*}
thus for $\gamma$ sufficiently close to 1, although $\|\theta_0-\Theta(0,\cdot)\|_1$ is small, $\|(v_0-V(0,\cdot),u_0-U(0,\cdot))\|_1$, $\|s_0(x)\|_1$ and $\|\chi_{0}-1\|_2$ can indeed be arbitrarily large. Therefore, in this sense Theorem 1.1 is a Nishida-Smoller-type result \cite{Nishida-Smoller-type} with large initial data.
\end{Remark}
\begin{Remark}
Compared with previous results on the large-time behavior of strong solutions to the 1D isentropic compressible NSAC system \cite{Chen-DCDS-B,Lei-Chen-JMAA,Lei-Chen-} in unbounded domain, Theorem 1.1 gives the first stability result of the composite wave to the Cauchy problem for the 1D compressible full NSAC system (\ref{1.9}) with arbitrarily large initial data. However, the premise for the validity of the above results is that $\gamma-1$ is sufficiently small. Therefore, it is still open when $\gamma$ is not close to $1$ on the global well-posedness and large-time behavior of solutions to the 1D isentropic/non-isentropic compressible NSAC system in unbounded domain so far, which is more difficult and challenging. These interesting problems will be pursued by the authors in the future.
\end{Remark}
\begin{Remark}
Motivated by the relationship between NSAC system and Navier-Stokes system, we assume that $\chi_{\pm}=1$. The case of $\chi_{+}\neq\chi_{-}$, which leads to more complex structure and difficulties,  is also left for future study.
\end{Remark}

\begin{Remark}
The Theorem 1.1 is also correct if constant coefficients $\nu$ and $\kappa$ are replaced by $\nu(\theta)=\nu\theta^\alpha$ and $\kappa(\theta)=\kappa\theta^\beta$, respectively, where $\alpha$ and $\beta$ are both positive constants.
\end{Remark}

Now we explain the main strategy of this paper. To prove Theorem 1.1, we mainly use the elementary $L^2$-energy method to derive the suitable a priori estimates of solutions to the initial value problem (\ref{2.1})-(\ref{2.2}). The key ingredient is to deduce the uniform-in-time positive lower and upper bounds on the specific volume $v(t,x)$ and temperature $\theta(t,x)$. Firstly, it is assumed that the parameter $\gamma-1$ is sufficiently small in this paper. This means that the uniform-in-time positive lower and upper bounds on the temperature $\theta(t,x)$ are easily accessible. Secondly, we apply the argument developed by Jiang \cite{Jiang-CMP} to derive a presentation of the specific volume $v(t,x)$ (see Lemma 2.5), which is obtained by using a special cut-off function when the far field states are different. Thanks to the choice of the cut-off function and delicate analysis based on the basic energy estimate, the uniform-in-time positive lower and upper bounds on the temperature $\theta(t,x)$, (\ref{2.38}) and the Gronwall inequality, the specific volume $v(t,x)$ is shown uniformly bounded from below and above with respect to the space and time (see Lemma 2.6). Finally, having obtained the uniformly positive lower and upper bounds of $v(t,x)$ and $\theta(t,x)$, the higher order energy estimates of solutions can be established by using the lower order energy estimates and Gronwall's inequality. Noticing that we need to deal with the high-order strongly nonlinear term $\left(\frac{\chi_x^2}{v^2}\right)_x$. Time independent estimate of $\|\phi_x(t)\|$ is a challenging endeavour, which only depends on the upper bound of $\|\chi_t(t)\|$ essentially (see Lemmas 2.8-2.9). Here the perturbed function $\phi(t,x):=v(t,x)-V(t,x)$. Then Theorem 2.1 can be proved by the standard continuation argument.

This paper is organized as follows. In Section 2, we first reformulate the system for the perturbation around the rarefaction waves in Section 1.1. Based on the a priori assumption (\ref{2.6-1}), we establish the a priori estimates of solutions to the system (\ref{2.1})-(\ref{2.2}), and then give the proof of Theorem 2.1 by the standard continuation argument.

\emph{\textbf{Notation.}} Throughout this paper, $C$ denotes some generic time-independent and $(\gamma-1)$-independent positive constant which may vary in different estimates. If the dependence needs to be explicitly pointed out, the notation $C(\cdot,\cdots,\cdot)$ or $C_i(\cdot,\cdots,\cdot)(i\in {\mathbb{N}})$ is used. For function spaces,  $L^p(\mathbb{R})(1\leq p\leq+\infty)$ represents the standard Lebesgue space with the norm  $\|\cdot\|_{L^p}:=\|\cdot\|_{L^p(\mathbb{R})}$,  and $W^{k,p}(\mathbb{R})(k\geq0)$ is the usual $k$-th order Sobolev space with its norm:
\[\|f\|_{k,p}=\left(\sum_{i=0}^{k}\left\|\partial_x^i f\right \|_{L^p}^2\right)^{\frac{1}{2}}.\]
For notational simplicity, when $p=2$, we denote $\|\cdot\|:=\|\cdot\|_{L^2}$ and $H^k(\mathbb{R}):=W^{k,2}(\mathbb{R})$ with its norm $\|\cdot\|_k:=\|\cdot\|_{k,2}$. Finally,  $\|\cdot\|_{L_{T,x}^\infty}$ stands for the norm $\|\cdot\|_{L^\infty([0,T]\times\mathbb{R})}$.

\section{Preliminaries}
\setcounter{equation}{0}

In this section, we give some basic properties of the rarefaction waves, which will be used later.

The unique weak entropy solution of Riemann problem $(\ref{1.25})$ is the rarefaction wave $w_i^{R}\left(\frac{x}{t}\right)$ defined by
\begin{eqnarray}\label{1.27}
w_i^{R}\left(\frac{x}{t}\right)=\left\{\begin{array}{ll}
w_{i-},\,\,x\leq w_{i-}t,\\[2mm]
\frac{x}{t},\,\,w_{i-}t\leq x\leq w_{i+}t,\\[2mm]
w_{i+},\,\,x\geq w_{i+}t.
\end{array}\right.
\end{eqnarray}

The above rarefaction wave $w_i^{R}\left(\frac{x}{t}\right)$ can be approximated by the unique global smooth solution $\omega_i(t,x)(i=1,3)$ of the following Cauchy problem:
\begin{eqnarray}\label{1.28}
\left\{\begin{array}{ll}
    \omega_{it}+\omega_i\omega_{ix}=0,~~t>0,~x\in\mathbb{R},\\[2mm]
    \omega_i(t,x)\big|_{t=0}=\omega_{i0}(x)=\frac{\omega_{i+}+\omega_{i-}}{2}+\frac{\omega_{i+}-\omega_{i-}}{2}K_q\int_0^{\varepsilon x}(1+y^2)^{-q}dy,
\end{array}\right.
\end{eqnarray}
where $\omega_{i+}$ and $\omega_{i-}(i=1,3)$ are given by $(\ref{1.26})$, $K_q>0$ is a positive constant, such that $K_q=\left(\int_0^{+\infty}(1+y^2)^{-q}dy\right)^{-1}$, $q>\frac{3}{2}$ and $\varepsilon>0$ are positive constants which will be specified later.

The properties of the solution $w_i(t,x)(i=1,3)$ of (\ref{1.28}) are listed in the following lemma.
\begin{Lemma}
Let each $i\in\{1,3\}$, $\widetilde{\omega}_i=\omega_{i+}-\omega_{i-}$, then the Cauchy problem (\ref{1.28}) admits a unique global smooth solution $\omega_i(t,x)$ satisfying:\\
(i) $\omega_{i-}<\omega_i(t,x)<\omega_{i+}$,~$\omega_{ix}(t,x)>0$ for any  $(t,x)\in\mathbb{R}^+\times\mathbb{R}$.\\
(ii) For any $p(1\leq p\leq+\infty)$, there exists a constant $C_{p,q}>0$ depending on $p,q$ such that for $t\geq 0$,
\begin{eqnarray*}
&&\|\omega_{ix}(t)\|_{L^p}\leq C_{p,q}\min\left\{\widetilde{\omega}_i\varepsilon^{1-\frac{1}{p}},\widetilde{\omega}_i^{\frac{1}{p}}t^{-1+\frac{1}{p}}\right\},\nonumber\\
&&\left\|\omega_{ixx}(t)\right\|_{L^p}\leq C_{p,q}\min\left\{\widetilde{\omega}_i\varepsilon^{2-\frac{1}{p}},\varepsilon^{\left(1-\frac{1}{2q}\right)\left(1-\frac{1}{p}\right)}\widetilde{\omega}_i^{-\frac{p-1}{2pq}}t^{-1-\frac{1}{2q}+\frac{1}{2pq}}\right\}.
\end{eqnarray*}
(iii) If $0<\omega_{i-}(<\omega_{i+})$ and $q$ is suitably large, then
\begin{equation*}
\left\{\begin{array}{ll}
     \left|\omega_i(t,x)-\omega_{i-}\right|\leq C\widetilde{\omega}_i\left(1+(\varepsilon x)^2\right)^{-\frac{q}{3}}\left(1+(\varepsilon \omega_{i-}t)^2\right)^{-\frac{q}{3}},~x\leq 0,\,\,t>0;\\[2mm]
     \left|\omega_{ix}(t,x)\right|\leq C\varepsilon\widetilde{\omega}_i\left(1+(\varepsilon x)^2\right)^{-\frac{q}{2}}\left(1+(\varepsilon \omega_{i+}t)^2\right)^{-\frac{q}{2}},~x\leq 0,\,\,t>0.
\end{array}\right.
\end{equation*}
(iv) If $(\omega_{i-}<)\omega_{i+}\leq 0$ and $q$ is suitably large, then
\begin{equation*}
\left\{\begin{array}{ll}
     \left|\omega_i(t,x)-\omega_{i+}\right|\leq C\widetilde{\omega}_i\left(1+(\varepsilon x)^2\right)^{-\frac{q}{3}}\left(1+(\varepsilon \omega_{i-}t)^2\right)^{-\frac{q}{3}},~x\geq 0,\,\,t>0;\\[2mm]
     \left|\omega_{ix}(t,x)\right|\leq C\varepsilon\widetilde{\omega}_i\left(1+(\varepsilon x)^2\right)^{-\frac{q}{2}}\left(1+(\varepsilon \omega_{i+}t)^2\right)^{-\frac{q}{2}},~x\geq 0,\,\,t>0.
\end{array}\right.
\end{equation*}
(v) \begin{equation*}
\lim_{t\rightarrow+\infty}\sup_{x\in\mathbb{R}}\left|\omega_i(t,x)-\omega_i^R\left(\frac{x}{t}\right)\right|=0.
\end{equation*}
\end{Lemma}

Based on Lemma 2.1 and the definition of  $(V,U,\Theta,S)(t,x)$ given in (\ref{1.29}) and (\ref{1.30}), we can deduce that

\begin{Lemma}
Let the strength of rarefaction waves $\delta:=|v_+-v_-|+|u_+-u_-|$, then the smooth approximate rarefaction wave $(V,U,\Theta,S)(t,x)$ constructed in (\ref{1.29})-(\ref{1.31}) has the following properties:\\
(i) $V_t(t,x)=U_x(t,x)>0,~\forall(t,x)\in\mathbb{R}^+\times\mathbb{R}.$\\
(ii) For any $p(1\leq p\leq+\infty)$, there exists a constant $C_{p,q}>0$ depending only on $p,q$ such that for $t\geq 0$,
\begin{eqnarray*}
&&\|(V_x,U_x,\Theta_x)(t)\|_{L^p}\leq C_{p,q}\min\left\{\delta\varepsilon^{1-\frac{1}{p}},\delta^{\frac{1}{p}}(t+t_0)^{-1+\frac{1}{p}}\right\},\nonumber\\
&&\left\|(V_{xx},U_{xx},\Theta_{xx})(t)\right\|_{L^p}\leq C_{p,q}\min\left\{\delta\varepsilon^{2-\frac{1}{p}},\delta^{-\frac{p-1}{2pq}}\varepsilon^{\left(1-\frac{1}{2q}\right)\left(1-\frac{1}{p}\right)}(t+t_0)^{-1-\frac{1}{2q}+\frac{1}{2pq}}\right\}.\\
\end{eqnarray*}

It is obvious that $\|V_x(t)\|^2_{L^2}$ is not integrable with respect to $t$, however we can get for any $r>0$ and $p>1$ that
\begin{eqnarray*}
&&\int_{0}^{\infty}\|(V_x,U_x,\Theta_x)(t)\|^{2+r}_{L^{2+r}}dt\leq C_{r}\delta,~~\forall r>0,\\
&&\int_{0}^{\infty}\left\|(V_{xx},U_{xx},\Theta_{xx})(t)\right\|_{L^p}dt\leq C_{p,q}\delta^{-\frac{p-1}{2pq}}\varepsilon^{\left(1-\frac{1}{2q}\right)\left(1-\frac{1}{p}\right)},~~\forall p>1.
\end{eqnarray*}
(iii) For each $p\geq 1$,
\begin{equation*}
\|(g(V,\Theta)_x,r(V,\Theta),q(V,\Theta))(t)\|_{L^p}\leq C(p,q)\delta^2\varepsilon^{1-\frac{1}{p}}(1+(\varepsilon t)^2+(\varepsilon t_0)^2)^{-\frac{q}{3}}.
\end{equation*}

Especially,
\begin{equation*}
\int_0^\infty\|(g(V,\Theta)_x,r(V,\Theta),q(V,\Theta))(t)\|_{L^p}dt\leq C(p)\delta^2\varepsilon^{2-\frac{1}{p}}.
\end{equation*}
(iv)\begin{equation*}
\lim_{t\rightarrow+\infty}\sup_{x\in\mathbb{R}}\left|\left(V,U,\Theta\right)(t,x)-\left(V^R,U^R,\Theta^R\right)\left(\frac{x}{t}\right)\right|=0.
\end{equation*}
(v)~$\left|\left(V_t,U_t,\Theta_t\right)(t,x)\right|\leq C\left|\left(V_x,U_x,\Theta_x\right)(t,x)\right|$.
\end{Lemma}

Then, according to (\ref{1.19}), (\ref{1.28}) and (\ref{1.29}), the approximate waves $(V,U,\Theta,S)(t,x)$ satisfies the following equations:
\begin{equation}\label{1.32}
\left\{\begin{array}{ll}
   V_t-U_x=0,\,\,\\[2mm]
   U_t+p(V,\Theta)_x=g(V,\Theta)_x,\,\,\\[2mm]
   \left[e(V,\Theta)+\frac{U^2}{2}\right]_t+(Up(V,\Theta))_x=q(V,\Theta),\,\,\\[2mm]
   C_v\Theta_t+p(V,\Theta)U_x=r(V,\Theta),\,\,\\[2mm]
   S_t(V,\Theta)=0,
\end{array}\right.
\end{equation}
where
\begin{eqnarray}\label{1.33}
g(V,\Theta)=&&p(V,\Theta)-p(V_1,\Theta_1)-p(V_3,\Theta_3)-p(v_m,u_m),\nonumber\\
q(V,\Theta)=&&\left(e(V,\Theta)-e(V_1,\Theta_1)-e(V_3,\Theta_3)\right)_t+\left(\frac{U^2}{2}-\frac{U_1^2}{2}-\frac{U_3^2}{2}\right)_t\\
&&+\left[Up(V,\Theta)-U_1p(V_1,\Theta_1)-U_3p(V_3,\Theta_3)\right]_x,\nonumber\\
r(V,\Theta)=&&p(V,\Theta)U_x-p(V_1,\Theta_1)U_{1x}-p(V_3,\Theta_3)U_{3x},\nonumber
\end{eqnarray}
and $\theta_m=\widetilde{\theta}(v_m,\overline{s})$.

\section{Reformulation of the original problem}
\setcounter{equation}{0}
In this section, we reformulate our original problem $(\ref{1.9})_1$-$(\ref{1.9})_4$ and (\ref{1.13}) into a perturbation one around the smooth approximation rarefaction wave $(V,U,\Theta,S)(t,x)$. An equivalent theorem of Theorem 1.1 (i.e., Theorem 3.1) is stated in this section. Then, based on the local existence result and the a priori estimates of the reformulated problem, the proof of Theorem 3.1 is presented at the end of Section 3.

 First of all, we define the perturbation functions $(\phi,\psi,\zeta,\varphi,\xi)(t,x)$ by
\begin{eqnarray*}
&&(\phi,\psi,\zeta,\varphi,\xi)(t,x)
=(v(t,x)-V(t,x),u(t,x)-U(t,x),\theta(t,x)-\Theta(t,x),s(t,x)-S(t,x),\chi(t,x)-1),
\end{eqnarray*}
then we have from $(\ref{1.9})_1$-$(\ref{1.9})_4$, (\ref{1.13}) and (\ref{1.32}) that
\begin{eqnarray}\label{2.1}
\left\{\begin{array}{ll}
    \phi_t-\psi_x=0,\,\,\\[2mm]
    \psi_t+\left[p(v,\theta)-p(V,\Theta)\right]_x=\nu\left(\frac{u_x}{v}\right)_x-\frac{1}{2}\left(\frac{\xi_x^2}{v^2}\right)_x-g(V,\Theta)_x,\,\,\\[2mm]
    \xi_t=-v\mu,\,\,\,\,\,\,~~~~~~~~~~~~~~~~~~~~~~~~~~~~~~~~~~~~~~~~~~~~~~~~~~~~~~~~~~`~~~~~~~~~~~~~\forall t>0, ~~x\in\mathbb{R}\\[2mm]
    \mu=-\left(\frac{\xi_x}{v}\right)_x+\xi(\xi+1)(\xi+2),\,\,\\[2mm]
    C_v\zeta_t+p(v,\theta)\psi_x+\left[p(v,\theta)-p(V,\Theta)\right]U_x=\kappa\left(\frac{\theta_x}{v}\right)_x+\nu\frac{u_x^2}{v}+v\mu^2-r(V,\Theta),\,\,\\[2mm]
    \varphi_t=\kappa\left(\frac{\theta_x}{v\theta}\right)_x+\kappa\frac{\theta_x^2}{v\theta^2}+\nu\frac{u_x^2}{v\theta}+\frac{v\mu^2}{\theta},
\end{array}\right.
\end{eqnarray}
with the initial condition:
\begin{eqnarray}\label{2.2}
&&(\phi,\psi,\zeta,\varphi,\xi)(t,x)\big|_{t=0}
=(\phi_0,\psi_0,\zeta_0,\varphi_0,\xi_0)(x)\\
:=&&(v_0(x)-V(0,x),u_0(x)-U(0,x),\theta_0(x)-\Theta(0,x),s_0(x)-S(0,x),\chi_0(x)-1),\,\,x\in\mathbb{R}.\nonumber
\end{eqnarray}

We define the following set of functions in which the solutions of the Cauchy problem (\ref{2.1})-(\ref{2.2}) is sought as follows:
\begin{eqnarray*}
&&X(0,T;M_0,M_1,N_0,N_1)
=\left\{(\phi,\psi,\zeta,\varphi,\xi)(t,x)\left|
\begin{array}{c}
(\phi,\psi,\zeta,\varphi,\xi)(t,x)\in C(0, T; \left(H^{2}(\mathbb{R})\right)^5),\\[2mm]
\phi_x(t,x)\in L^2(0, T; H^{1}(\mathbb{R})),\\[2mm]
(\psi_x,\zeta_x)(t,x)\in L^2(0, T; \left(H^{2}(\mathbb{R})\right)^2),\\[2mm]
\xi(t,x)\in L^2(0, T; H^{3}(\mathbb{R})),\\[2mm]
 \displaystyle M_0\leq v(t,x)\leq M_1,~N_0\leq \theta(t,x)\leq N_1,
\end{array}
\right.\right\}
 \end{eqnarray*}
where $T>0$ is a positive constant. $M_0,M_1,N_0,N_1$ and $T$ are some positive constants.

Then to prove Theorem 1.1, it suffices to show the following theorem.
\begin{Theorem}
Under the assumptions of Theorem 1.1, there exists a sufficiently small positive constant $\varepsilon_1> 0$ depending only on $\overline{V},~\underline{V},~\underline{\Theta},~\overline{\Theta}$, $\delta$,~$\|\phi_0\|_2,~\|\psi_0\|_2$, $\|\frac{\zeta_0}{\sqrt{\gamma-1}}\|_2$ and $\|\xi_0\|_2$, such that if $0< \varepsilon_0\leq\varepsilon_1$, then the Cauchy problem (\ref{2.1})-(\ref{2.2}) admits a unique global-in-time smooth solution $(\phi,\psi,\zeta,\varphi,\xi)(t,x)$, which satisfies that
\begin{equation}\label{2.3}
C_0^{-1}\leq v(t,x)\leq C_0,~\underline{\Theta}\leq \theta(t,x)\leq\overline{\Theta},~0\leq\chi(t,x)\leq 1,~\forall(t,x)\in[0,+\infty)\times\mathbb{R},
\end{equation}
\begin{eqnarray}\label{2.4}
&&\left\|\left(\phi,\psi,\frac{\zeta}{\sqrt{\gamma-1}},\varphi,\xi\right)(t)\right\|_2^2+\int_0^t\left(\|\varphi_x(\tau)\|_1^2+\|(\psi_x,\zeta_x)(\tau)\|_2^2+\|\xi(\tau)\|_3^2\right)d\tau\nonumber\\
\leq&&C_2\left(\left\|\left(\phi_0,\psi_0,\frac{\zeta_0}{\sqrt{\gamma-1}},\xi_0\right)\right\|_2^2+1\right),~\forall t>0.
\end{eqnarray}
Here $C_0$ is a positive constant depending only on $\overline{V},~\underline{V},~\underline{\Theta},~\overline{\Theta}$,~$\delta$,~$\|\phi_0\|,~\|\psi_0\|$, $\|\frac{\zeta_0}{\sqrt{\gamma-1}}\|$ and $\|\xi_0\|_1$, and $C_2$ is a positive constant depending only on $\overline{V},~\underline{V},~\underline{\Theta},~\overline{\Theta}$,~$\delta$,~$\|\phi_0\|_2,~\|\psi_0\|_2$, $\left\|\frac{\zeta_0}{\sqrt{\gamma-1}}\right\|_2$ and $\|\xi_0\|_2$.

Furthermore, the following large-time behaviors hold:
\begin{equation}\label{2.5}
\lim_{t\rightarrow+\infty}\|(\phi,\psi,\zeta,\varphi,\xi)(t)\|_{L^{\infty}(\mathbb{R})}=0.
\end{equation}
\end{Theorem}

Indeed, once Theorem 3.1 is proved, (\ref{1.36}) follows immediately from (\ref{2.5}) and Lemma 2.2 $(iv)$. The other conclusions of Theorem 1.1 follow directly from Theorem 3.1.

Theorem 3.1 will be obtained by the standard  continuation argument, which is based on the following local existence result Proposition 3.1 and the a priori estimates stated in Proposition 3.2.

\begin{Proposition} [Local existence]
Suppose that  the conditions  of  Theorem 1.1 hold, and $\nu_0\leq v_0(x)\leq\nu_1$ and $\varsigma_0\leq \theta_0(x)\leq\varsigma_1$ for some positive constants $\nu_0,\nu_2,\varsigma_0,\varsigma_1>0$ and all $x\in\mathbb{R}$, then there exists a small positive constant $t_0=t_0(\nu_0,\varsigma_0, N_0)$ depending only on $\nu_1$ and $\biggl\|\biggl(v_0-V(\cdot,0),u_0-U(\cdot,0),
\frac{\theta_0-\Theta(\cdot,0)}{\sqrt{\gamma-1}},\chi_0-1\biggr)\biggr\|_2$, such that the Cauchy problem (\ref{2.1})-(\ref{2.2}) admits a unique smooth solution $(\phi,\psi,\zeta,\varphi,\xi)(t,x)\in X\left(0,t_0; \frac{1}{2}\nu_0, 2\nu_1,\frac{1}{2}\varsigma_0,2\varsigma_1\right)$, which satisfies that
\begin{eqnarray*}
0<\underline{V}\leq v(t,x)\leq \overline{V},~\underline{\Theta}\leq \theta(t,x)\leq \overline{\Theta},~ 0\leq\chi(t,x)\leq1,~
\forall(t,x)\in[0,t_0]\times\mathbb{R},
\end{eqnarray*}
and
\begin{eqnarray*}
\sup_{0\leq t\leq t_0}\left\|\left(\phi,\psi,\frac{\zeta}{\sqrt{\gamma-1}},\varphi,\xi\right)(t)\right\|_{2}^{2}
 \leq\vartheta\left\|\left(\phi_0,\psi_0,\frac{\zeta_0}{\sqrt{\gamma-1}},\xi_0\right)\right\|_{2}^{2}.
\end{eqnarray*}
where $\vartheta=\vartheta(\nu_0,\nu_1,\varsigma_0,\varsigma_1)>1$ is a constant depending only on $\nu_0$, $\nu_1$, $\varsigma_0$ and $\varsigma_1$.
\end{Proposition}

\begin{Proposition} [A priori estimates]
Under the conditions of  Theorem 3.1 hold, suppose that $(\phi,\psi,\zeta,\varphi,\xi)(t,x)\in X(0,T;M_0,M_1,N_0,N_1)$ is a solution of the Cauchy problem (\ref{2.1})-(\ref{2.2}) in the strip $\Pi_T=[0,T]\times\mathbb{R}$ for some positive constants $M_0$, $M_1$, $N_0$, $N_1$ and $T$.
Assume that the following a priori assumptions
\begin{eqnarray}\label{2.6-1}
      &&0<M_0\leq v(t,x)\leq M_1,~~\underline{\Theta}\leq\theta(t,x)\leq\overline{\Theta},~~\\ &&\left\|\left(\phi,\psi,\frac{\zeta}{\sqrt{\gamma-1}},\xi\right)(t)\right\|_2^2+\int_0^t\left(\|\phi_x(\tau)\|_1^2+\|(\psi_x,\zeta_x)(\tau)\|_2^2+\|\xi(\tau)\|_3^2\right)d\tau\leq N^2,\nonumber
\end{eqnarray}
hold for all $x\in\mathbb{R}$ and $ 0\leq t\leq T$, and some generic positive constants $M_0,M_1,N$. The definitions of $\underline{\Theta}$ and $\overline{\Theta}$ have already been given in the conditions of Theorem 1.1.

Then, there exist a positive constant $C_0$ depending only on $\underline{V}, \overline{V}$, $\underline{\Theta}$, $\overline{\Theta}$, $\delta$, $\|\phi_0\|$, $\|\psi_0\|$, $\left\|\frac{\zeta_0}{\sqrt{\gamma-1}}\right\|$ and $\|\xi_0\|_1$, and a positive constant $C_3$ depending only on $\underline{V}, \overline{V}$, $\underline{\Theta}$, $\overline{\Theta}$, $\delta$, $\|\phi_0\|_2$, $\|\psi_0\|_2$, $\left\|\frac{\zeta_0}{\sqrt{\gamma-1}}\right\|_{2}$ and $\|\xi_0\|_2$, such that the following two estimates hold:
\begin{eqnarray}
&&C_0^{-1}\leq v(t,x)\leq C_0,~\underline{\Theta}\leq \theta(t,x)\leq \overline{\Theta},\quad 0\leq\chi(t,x)\leq1,\,\, \forall\,(t,x)\in[0,T]\times\mathbb{R},\label{2.6}\\
&&\left\|\left(\phi,\psi,\frac{\zeta}{\sqrt{\gamma-1}},\varphi,\xi\right)(t)\right\|^2_{2}+\int_0^t\left(\|\phi_x(\tau)\|^2_{1}+\|(\psi_x,\zeta_x)(\tau)\|^2_{2}+\|\xi(\tau)\|_3^2\right)d\tau\nonumber\\
\leq&&C_3\left(\left\|\left(\phi_0,\psi_0,\frac{\zeta_0}{\sqrt{\gamma-1}},\xi_0\right)\right\|^2_{2}+1\right), \,\, \forall\, t\in[0,T].\label{2.7}
\end{eqnarray}
\end{Proposition}

The proof of Proposition 3.1 can be done by using the standard iteration technique, which is similar to that of Theorem 3.1 in \cite{Ding-2013}. The details are omitted here for brevity. Hence to prove Theorem 3.1, it remains to show Proposition 3.2, which is left in the following subsection.
Finally, we present the proof of  Theorem 2.1 as follows.

Now we prove the Theorem 3.1 by using the local existence Proposition 3.1 and the a priori
estimate Proposition 3.2.\\
\noindent{\bf Proof of Theorem 3.1.} From Proposition 3.2, we derive
\begin{equation*}
|\zeta(t,x)|\leq C\sqrt{\gamma-1}\left\|\frac{\zeta}{\sqrt{\gamma-1}}\right\|^{\frac{1}{2}}\left\|\frac{\zeta_x}{\sqrt{\gamma-1}}\right\|^{\frac{1}{2}}\leq C_3\left(\left\|\left(\phi_0,\psi_0,\frac{\zeta_0}{\sqrt{\gamma-1}},\xi_0\right)\right\|_{2}+1\right)\sqrt{\gamma-1}.
\end{equation*}
Since
\begin{equation*}
\theta(t,x)=\Theta(t,x)+\zeta(t,x),~~2\underline{\Theta}\leq\Theta(t,x)\leq\frac{1}{2}\overline{\Theta},
\end{equation*}
then letting
\begin{eqnarray}\label{2.65}
C_3\left(\left\|\left(\phi_0,\psi_0,\frac{\zeta_0}{\sqrt{\gamma-1}},\xi_0\right)\right\|_{2}+1\right)\sqrt{\gamma-1}\leq\min\{\underline{\Theta},\frac{1}{2}\overline{\Theta}\}, \end{eqnarray}
we obtain
\begin{equation*}
\underline{\Theta}\leq\theta(t,x)\leq\overline{\Theta}.
\end{equation*}
Furthermore, according to Proposition 3.2, we can check that
\begin{equation*}
\int_0^{+\infty}\left(\|(\phi_x,\psi_x,\zeta_x,\xi_x)(t)\|^2+\left|\frac{d}{dt}\|(\phi_x,\psi_x,\zeta_x,\xi_x)(t)\|^2\right|\right)dx\leq +\infty,
\end{equation*}
which means
\begin{equation*}
\|(\phi_x,\psi_x,\zeta_x,\xi_x)(t,\cdot)\|\rightarrow 0, \textnormal{as}~t\rightarrow\infty.
\end{equation*}
Consequently, we have
\begin{equation*}
\|(\phi,\psi,\zeta,\xi)(t,\cdot)\|_{L^\infty}\leq C \|(\phi,\psi,\zeta,\xi)(t,\cdot)\|^\frac{1}{2}\|(\phi_x,\psi_x,\zeta_x,\xi_x)(t,\cdot)\|^\frac{1}{2}\rightarrow 0,~\textnormal{as}~t\rightarrow\infty.
\end{equation*}
This complete the proof of the Theorem 3.1.

\section{Proof of the a priori estimates}
\setcounter{equation}{0}

The main target of this section is to deduce the a priori estimates (\ref{2.6})-(\ref{2.7}) for the Cauchy problem (\ref{2.1})-(\ref{2.2}).

Before proving the a priori estimates, we first list the assumptions on the pressure function $p(v,\theta)$ and the internal energy $e(v,\theta)$ used throughout this paper:
\begin{equation}\label{2.8-1}
p_v(v,\theta)=\frac{\partial p(v,\theta)}{\partial v}<0,~~~e_\theta(v,\theta)=\frac{\partial e(v,\theta)}{\partial\theta}>0,
\end{equation}
and
\begin{equation}\label{2.8-2}
\widetilde{p}_{vv}(v,s)=\frac{\partial^2 \widetilde{p}(v,s)}{\partial v^2}>0~\textnormal{and}~\widetilde{p}(v,s)~\textnormal{is~convex~with~respect~to}~(v,s).
\end{equation}

From (\ref{1.12}) and (\ref{2.8-1}), we can deduce that
\begin{eqnarray}\label{2.8-3}
\widetilde{p}_v(v,s)=p_v(v,\theta)-\frac{\theta (p_\theta(v,\theta))^2}{e_\theta(v,\theta)}<0,
\end{eqnarray}
\begin{eqnarray*}
\left\{\begin{array}{ll}
     \widetilde{e}_{ss}(v,s)=\frac{\theta}{e_\theta(v,\theta)}>0,\\[2mm]
     \widetilde{e}_{vs}(v,s)=\frac{\theta p_\theta(v,\theta)}{e_\theta(v,\theta)},\\[2mm]
     \widetilde{e}_{vv}(v,s)=-p_v(v,\theta)+\frac{\theta (p_\theta(v,\theta))^2}{e_\theta(v,\theta)}>0,
\end{array}
\right.
 \end{eqnarray*}
and
\begin{eqnarray}\label{2.8-4}
\widetilde{e}_{ss}(v,s)\widetilde{e}_{vv}(v,s)-(\widetilde{e}_{vs}(v,s))^2=-\frac{\theta p_v(v,\theta)}{e_\theta(v,\theta)}>0.
\end{eqnarray}

From (\ref{2.8-4}), we know that $\widetilde{e}(v,s)$ is convex with respect to $v$ and $s$. Consequently, $\widetilde{e}(v,s)+\frac{u^2}{2}$ is a strictly convex function of $(v,u,s)$. Now we can construct the following normalized entropy $\eta(v,u,s;V,U,S)$ around $(V,U,S)(t,x)$, which is the smooth approximation of the rarefaction waves $(V^R,U^R,S^R)\left(\frac{x}{t}\right)$ which consists of the 1-rarefaction wave $\left(V_1^R\left(\frac{x}{t}\right),U_1^R\left(\frac{x}{t}\right),\overline{s}\right)$ and 3-rarefaction wave $\left(V_3^R\left(\frac{x}{t}\right),U_3^R\left(\frac{x}{t}\right),\overline{s}\right)$:
\begin{eqnarray}\label{2.8-5}
\eta(v,u,s;V,U,S)=&&\left(e(v,\theta)+\frac{u^2}{2}\right)-\left(e(V,\Theta)+\frac{U^2}{2}\right)\nonumber\\
&&-\left\{-p(V,\Theta)(v-V)+U(u-U)+\Theta(s-S)\right\},
\end{eqnarray}
here we have used the fact that $\widetilde{e}_{v}(v,s)=-p(v,\theta),\widetilde{e}_{s}(v,s)=\theta$. Note that for the ideal polytropic gas, $p(v,\theta)$ and $e(v,\theta)$ satisfy the constitutive relations $(\ref{1.15})$.

Next, we prove the a priori estimates (\ref{2.6})-(\ref{2.7}) via a series of lemmas.  First, we give the following key lemma.
\begin{Lemma}(Basic energy estimates) Under the assumptions of Proposition 3.2, then there exists a uniform positive constant $C_4>0$ depending only on $\underline{V},\overline{V}$, $\underline{\Theta}$, $\overline{\Theta}$ and $\delta$, it holds for $t\in[0,T]$ that
\begin{eqnarray}\label{2.8}
     &&\int_{\mathbb{R}}\left(R\Theta\Phi\left(\frac{v}{V}\right)+\frac{\psi^2}{2}+C_v\Theta\Phi\left(\frac{\theta}{\Theta}\right)+\frac{\xi_x^2}{2v}+\frac{\xi^2(\xi+2)^2}{4}\right)dx\nonumber\\
     &&+\int_0^t\int_{\mathbb{R}}\left(\frac{v\Theta}{\theta}\mu^2+\nu\frac{\Theta}{v\theta}\psi_x^2+\kappa\frac{\Theta}{v\theta^2}\zeta_x^2\right)dxd\tau\nonumber\\
     &&+\int_0^t\int_{\mathbb{R}}\left[\widetilde{p}(v,s)-\widetilde{p}(V,\overline{s})-\widetilde{p}_v(V,\overline{s})\phi-\widetilde{p}_s(V,\overline{s})\varphi\right]U_xdxd\tau\\
     \leq&&C_4\left(\left\|\left(\phi_0,\psi_0,\frac{\zeta_0}{\sqrt{\gamma-1}}\right)\right\|^2+\|\xi_{0}\|_1^2+1\right),\nonumber
\end{eqnarray}
provided that the parameters $\varepsilon>0$ and $\gamma-1$ are sufficiently small. Here the entropy $\eta(v,u,s;V,U,S)$ defined by $(\ref{2.8-5})$ takes the form
\begin{eqnarray}\label{2.8-7}
\eta(v,u,\theta;V,U,\Theta)=&&R\Theta\Phi\left(\frac{v}{V}\right)+\frac{1}{2}(u-U)^2+\frac{R}{\gamma-1}\Theta\Phi\left(\frac{\theta}{\Theta}\right),
\end{eqnarray}
where the function $\Phi(s)$ is defined by
\begin{equation}\label{2.8-8}
\Phi(s)=s-1-\ln s.
\end{equation}
\end{Lemma}
\noindent{\bf Proof.}~~From $(\ref{1.15})$, $(\ref{2.8-7})$ and $(\ref{2.8-8})$, we can get
\begin{eqnarray}\label{2.9}
&&\eta_t+\nu\frac{\Theta}{v\theta}\psi_x^2+\kappa\frac{\Theta}{v\theta^2}\zeta_x^2+\left[\widetilde{p}(v,s)-\widetilde{p}(V,\overline{s})-\widetilde{p}_v(V,\overline{s})\phi-\widetilde{p}_s(V,\overline{s})\varphi\right]U_x\nonumber\\
=&&\left\{\nu\frac{\psi\psi_x}{v}+\kappa\frac{\zeta\zeta_x}{v\theta}-\psi[p(v,\theta)-p(V,\Theta)]\right\}_x+\underbrace{\left(-\nu\frac{\psi U_x\phi_x}{v^2}+2\nu\frac{\zeta U_x\psi_x}{v\theta}-\kappa\frac{\Theta_x\zeta\phi_x}{v^2\theta}+\kappa\frac{\Theta_x\zeta\zeta_x}{v\theta^2}\right)}_{\mathcal{R}_1}\nonumber\\
&&+\underbrace{\left(\nu\frac{U_{xx}\psi}{v}+\kappa\frac{\zeta\Theta_{xx}}{v\theta}\right)}_{\mathcal{R}_2}+\underbrace{\left(-\nu\frac{V_xU_x\psi}{v^2}+\nu\frac{\zeta U_x^2}{v\theta}-\kappa\frac{\zeta\Theta_xV_x}{v^2\theta}\right)}_{\mathcal{R}_3}\\
&&-\underbrace{\frac{1}{C_v}r(V,\Theta)\varphi-q(V,\Theta)+Ug(V,\Theta)_x-\psi g(V,\Theta)_x}_{\mathcal{R}_4}-\frac{1}{2}\left(\frac{\xi_x^2}{v^2}\right)_x\psi+\frac{v}{\theta}\zeta\mu^2.\nonumber
\end{eqnarray}
We multiply $(\ref{2.1})_3$ by $\mu$ to obtain
\begin{eqnarray}\label{2.10}
\xi_t\mu=-v\mu^2.
\end{eqnarray}
A direct calculation yields
\begin{displaymath}
    \begin{split}
    \xi_t\mu=&\xi_t\left[-\left(\frac{\xi_x}{v}\right)_x+\xi(\xi+1)(\xi+2)\right]\\
    =&-\left(\frac{\xi_t\xi_x}{v}\right)_x+\left(\frac{\xi_x^2}{2v}+\frac{\xi^2(\xi+2)^2}{4}\right)_t+\frac{\xi_x^2U_x}{2v^2}+\frac{\xi_x^2\psi_x}{2v^2},
 \end{split}
\end{displaymath}
that is,
\begin{eqnarray}\label{2.11}
\left[\frac{\xi_x^2}{2v}+\frac{\xi^2(\xi+2)^2}{4}\right]_t+v\mu^2=\left(\frac{\xi_t\xi_x}{v}\right)_x-\frac{\xi_x^2U_x}{2v^2}-\frac{\xi_x^2\psi_x}{2v^2}.
\end{eqnarray}

Combining (\ref{2.9}) with (\ref{2.11}), then integrating the resultant equation over $(t,x)\in[0,t]\times\mathbb{R}$ yields that
\begin{eqnarray}\label{2.12}
&&\int_\mathbb{R}\left(\eta+\frac{\xi_x^2}{2v}+\frac{\xi^2(\xi+2)^2}{4}\right) dx+\int_0^t\int_\mathbb{R}\left(\nu\frac{\Theta}{v\theta}\psi_x^2+\kappa\frac{\Theta}{v\theta^2}\zeta_x^2+\frac{v\Theta}{\theta}\mu^2\right)dxd\tau\nonumber\\
&&+\int_0^t\int_\mathbb{R}\left[\widetilde{p}(v,s)-\widetilde{p}(V,\overline{s})-\widetilde{p}_v(V,\overline{s})\phi-\widetilde{p}_s(V,\overline{s})\varphi\right]U_xdxd\tau+\int_0^t\int_\mathbb{R}\frac{\xi_x^2U_x}{2v^2}dxd\tau\nonumber\\
=&&\int_\mathbb{R}\left(\eta_0+\frac{\xi_{0x}^2}{2v}+\frac{\xi_0^2(\xi_0+2)^2}{4}\right)dx+\sum_{i=1}^4\int_0^t\int_\mathbb{R}\mathcal{R}_idxd\tau.
\end{eqnarray}

By applying the H\"{o}lder inequality, the Cauchy inequality, the Young inequality, $\gamma-1$ sufficiently small and Lemma 2.2, we can get
\begin{eqnarray}\label{2.13}
\int_0^t\int_\mathbb{R}\mathcal{R}_1dxd\tau\leq&&\frac{1}{8}\int_0^t\int_\mathbb{R}\left(\nu\frac{\Theta}{v\theta}\psi_x^2+\kappa\frac{\Theta}{v\theta^2}\zeta_x^2\right)dxd\tau+C\int_0^t\int_\mathbb{R}\left(\frac{\zeta^2U_x^2}{v\theta}+\frac{\Theta_x^2\zeta^2}{v\theta^2}\right)dxd\tau\nonumber\\
&&+C(M_0^{-1},\underline{\Theta})\int_0^t\int_\mathbb{R}(|\psi\phi_xU_x|+|\zeta\Theta_x\phi_x|)dxd\tau\\
\leq&&\frac{1}{8}\int_0^t\int_\mathbb{R}\left(\nu\frac{\Theta}{v\theta}\psi_x^2+\kappa\frac{\Theta}{v\theta^2}\zeta_x^2\right)dxd\tau+\varepsilon^{\frac{1}{4}}\int_0^t\|\phi_x(\tau)\|^2d\tau\nonumber\\
&&+C_5(M_0^{-1},\underline{\Theta})\varepsilon^{\frac{1}{4}}\int_0^t\left\|\left(\psi,\frac{\zeta}{\sqrt{\gamma-1}}\right)(\tau)\right\|^2(1+\tau)^{-\frac{7}{4}}d\tau,\nonumber
\end{eqnarray}
\begin{eqnarray}\label{2.14}
\int_0^t\int_\mathbb{R}\mathcal{R}_2dxd\tau
\leq&&C\int_0^t\left(\|U_{xx}(\tau)\|\left\|\frac{\psi}{v}(\tau)\right\|+\|\Theta_{xx}(\tau)\|\left\|\frac{\zeta}{v\theta}(\tau)\right\|\right)d\tau\\
\leq&&C\int_0^t\left(\|U_{xx}(\tau)\|+\|\Theta_{xx}(\tau)\|\right)d\tau+C(M_0^{-1},\underline{\Theta})\int_0^t\|(U_{xx},\Theta_{xx})(\tau)\|\|(\psi,\zeta)(\tau)\|^2d\tau\nonumber\\
\leq&&C\delta^{-\frac{1}{4q}}\varepsilon^{\frac{2q-1}{4q}}+C_6(M_0^{-1},\underline{\Theta})\varepsilon^{\frac{2q-1}{4q}}\int_0^t\left\|\left(\psi,\frac{\zeta}{\sqrt{\gamma-1}}\right)(\tau)\right\|^2(1+\tau)^{-1-\frac{1}{4q}}d\tau,\nonumber
\end{eqnarray}
\begin{eqnarray}\label{2.15}
\int_0^t\int_\mathbb{R}\mathcal{R}_3dxd\tau
\leq&&C\int_0^t\int_\mathbb{R}\left(|V_x|^{\frac{5}{2}}+|U_x|^{\frac{5}{2}}\right)dxd\tau+C\int_0^t\int_\mathbb{R}\left(\left|\frac{U_x\psi}{v^2}\right|+\left|\frac{U_x\zeta}{v\theta}\right|+\left|\frac{\Theta_x\zeta}{v^2\theta}\right|\right)dxd\tau\nonumber\\
\leq&&C\int_0^t\|(V_{x},U_x,\Theta_{x})(\tau)\|^{\frac{5}{2}}_{L^{\frac{5}{2}}}d\tau+C(M_0^{-1},\underline{\Theta})\int_0^t|(U_{x},\Theta_{x})|^{\frac{3}{2}}|(\psi,\zeta)|^2d\tau\\
\leq&&C\delta+C_7(M_0^{-1},\underline{\Theta})\varepsilon^{\frac{1}{6}}\int_0^t\left\|\left(\psi,\frac{\zeta}{\sqrt{\gamma-1}}\right)(\tau)\right\|^2(1+\tau)^{-\frac{4}{3}}d\tau,\nonumber
\end{eqnarray}
and
\begin{eqnarray}\label{2.16}
&&\int_0^t\int_\mathbb{R}\mathcal{R}_4dxd\tau\nonumber\\
\leq&&C(\underline{V},\overline{V},\underline{\Theta},\overline{\Theta})\int_0^t\left(\|r(V,\Theta)(\tau)\|\|\varphi(\tau)\|+\|(q(V,\Theta),g(V,\Theta)_x)(\tau)\|_{L^1}+\|g(V,\Theta)_x(\tau)\|\|\psi(\tau)\|\right)d\tau\nonumber\\
\leq&&C(\underline{V},\overline{V},\underline{\Theta},\overline{\Theta})\int_0^t\|r(V,\Theta)(\tau)\|\|\varphi(\tau)\|^2d\tau+C\int_0^t\|g(V,\Theta)_x(\tau)\|\|\psi(\tau)\|^2d\tau\\
&&+C(\underline{V},\overline{V},\underline{\Theta},\overline{\Theta})\int_0^t\left(\|(r(V,\Theta),g(V,\Theta)_x)(\tau)\|+\|(q(V,\Theta),g(V,\Theta)_x)(\tau)\|_{L^1}\right)d\tau\nonumber\\
\leq&&C(\underline{V},\overline{V},\underline{\Theta},\overline{\Theta})\int_0^t\|(\varphi,\psi)(\tau)\|^2(1+\tau)^{-\frac{3}{2}}d\tau+C(\underline{V},\overline{V},\underline{\Theta},\overline{\Theta})\delta^2.\nonumber
\end{eqnarray}

If we choose $\varepsilon$ sufficiently small such that
\begin{eqnarray}\label{2.16-1}
   &&C_5(M_0^{-1},\underline{\Theta})\varepsilon^{\frac{1}{4}}\leq 1,~C_6(M_0^{-1},\underline{\Theta})\varepsilon^{\frac{2q-1}{4q}}\leq 1,~~\varepsilon^{\frac{1}{4}}N^2\leq 1,~C_7(M_0^{-1},\underline{\Theta})\varepsilon^{\frac{1}{6}}\leq 1,
\end{eqnarray}
then we obtain (\ref{2.8}) by inserting (\ref{2.13})-(\ref{2.16}) into (\ref{2.12}), exploiting the Gronwall inequality and $(\ref{2.6-1})_3$. This complete the proof of Lemma 4.1.

In Lemma 4.1, the smallness of $\varepsilon$ and $\gamma-1$ is used to guarantee that the basic energy (\ref{2.8}) is  bounded only by the initial data. Based on the basic energy estimate, we shall show that the specific volume $v(t,x)$ is uniformly bounded from below and above, which in turn actually decides the smallness of $\varepsilon$ and $\gamma-1$. This is why we need such the initial condition. To prove Proposition 3.2, we first aim to derive the uniform bounds of $v(t,x)$. We introduce two new variables: $\widetilde{v}=\frac{v}{V}$ and $\widetilde{\theta}=\frac{\theta}{\Theta}$. Then we have the following result.
\begin{Lemma}
Let $\alpha_1,\alpha_2$ be the two positive roots of the equation $y-\ln y-1=C_4$, and the constant $C_4$ be the same in (\ref{2.8}). Then,
\begin{equation}\label{2.18}
\alpha_1\leq \int_i^{i+1}\widetilde{v}(t,x)dx,~~\int_i^{i+1}\widetilde{\theta}(t,x)dx\leq \alpha_2,~~t\geq 0,
\end{equation}
and for each $t\geq 0$, there are points $a_i(t),b_i(t)\in[i,i+1]$ such that
\begin{equation}\label{2.19}
\alpha_1\leq \widetilde{v}(t,a_i(t)),~\widetilde{\theta}(t,b_i(t))\leq \alpha_2,~~t\geq 0,
\end{equation}
where $i=0,\pm 1,\pm 2,\cdots$.
\end{Lemma}
\noindent{\bf Proof.}~~(\ref{2.18}) and (\ref{2.19}) are similar to Lemma 3.3 in \cite{Huang-Wang-2016}. Thus we omit here for brevity. This complete the proof of Lemma 4.2.

The following lemma can be found in Jiang [\cite{Jiang-CMP}, Lemma 2.3].
\begin{Lemma}
For each $x\in[i,i+1],i=0,\pm 1,\pm 2,\cdots$, it follows from $(\ref{1.9})_2$ that
\begin{eqnarray}\label{2.20}
v(t,x)=B(t,x)Y(t)+\frac{1}{\nu}\int_0^t\frac{B(t,x)Y(t)}{B(\tau,x)Y(\tau)}\left(R\theta+\frac{\chi_x^2}{2v}\right)(\tau,x)d\tau,
\end{eqnarray}
where
\begin{eqnarray}\label{2.21}
&&B(t,x)=v_0(x)\exp\left\{\int_x^{+\infty}\frac{1}{\nu}(u_0(y)-u(t,y))\beta(y)dy\right\},\nonumber\\
&&Y(t)=\exp\left\{\frac{1}{\nu}\int_0^t\int_{i+1}^{i+2}\omega(\tau,x)dxd\tau\right\},\\
&&A(t,x)=\exp\left\{\frac{1}{\nu}\int_0^t\left(R\frac{\theta}{v}+\frac{\chi_x^2}{2v^2}\right)(\tau,x)d\tau\right\},\nonumber\\
&&\mathcal{E}(t,x)=\nu\frac{u_x}{v}-\frac{\chi_x^2}{2v^2}-\frac{R\theta}{v},\nonumber
\end{eqnarray}
and let $\beta(x)\in W^{1,\infty}(\mathbb{R})$ be defined by
\begin{equation}\label{2.22}
\beta(x)=\left\{\begin{array}{ll}
     1,~~~~~~~~~~~~~x\leq i+1,\\[2mm]
     i+2-x,~~~~i+1\leq x\leq i+2,\\[2mm]
     0,~~~~~~~~~~~~~x\geq i+2.
     \end{array}\right.
\end{equation}
\end{Lemma}
\noindent{\bf Proof.}~~We multiply $(\ref{1.9})_2$ by $\beta(x)$ defined in $(\ref{2.22})$ to obtain
\begin{equation*}
(\beta u)_t=(\mathcal{E}\beta)_x-\beta_x\mathcal{E}.
\end{equation*}

Integrating the above equation over $(x,\infty)(x\in[i,i+1])$ with respect to $x$, recalling $(\ref{1.9})_1$ and the definition of $\beta(x)$ and $\mathcal{E}(t,x)$, we can get
\begin{eqnarray}\label{2.23}
\int_x^{+\infty}(\beta u)_tdy=&&-\mathcal{E}-\int_x^{+\infty}\mathcal{E}\beta_xdy\nonumber\\
=&&-\nu(\ln v)_t+R\frac{\theta}{v}+\frac{\chi_x^2}{2v^2}+\int_{i+1}^{i+2}\mathcal{E}(t,y)dy.
\end{eqnarray}

Then integrating (\ref{2.23}) over $(0,t)$ with respect to $t$ and then taking the exponential on both sides of resultant equation yield that
\begin{eqnarray}\label{2.24}
v(t,x)=B(t,x)Y(t)\exp\left\{\frac{1}{\nu}\int_0^t\left(R\frac{\theta}{v}+\frac{\chi_x^2}{2v^2}\right)(\tau,x)d\tau\right\},
\end{eqnarray}
where $B(t,x)$ and $Y(t)$ are given in $(\ref{2.21})_1$ and $(\ref{2.21})_2$, respectively.

Note that
\begin{eqnarray}\label{2.25}
\frac{d}{dt}A(t,x)=\frac{A(t,x)}{\nu}\left(R\frac{\theta}{v}+\frac{\chi_x^2}{2v^2}\right)(t,x)=\frac{1}{\nu B(t,x)Y(t)}\left(R\theta+\frac{\chi_x^2}{2v}\right)(t,x),
\end{eqnarray}
where $A(t,x)$ is given by $(\ref{2.21})_3$. Later we integrate $(\ref{2.25})$ over $(0,t)$ with respect to $t$, which yield that
 \begin{eqnarray}\label{2.26}
A(t,x)=1+\int_0^t\frac{1}{\nu B(\tau,x)Y(\tau)}\left(R\theta+\frac{\chi_x^2}{2v}\right)(\tau,x)d\tau.
\end{eqnarray}

We get immediately $(\ref{2.20})$ by inserting $(\ref{2.26})$ into $(\ref{2.24})$. This complete the proof of Lemma 4.3.\\
\begin{Lemma} Under the assumptions of Proposition 3.2, then there exists a uniform positive constant $C_0>0$ depending only on $\underline{V},\overline{V}$, $\underline{\Theta}$, $\overline{\Theta}$, $\delta$, $\|\phi_{0}\|$, $\|\psi_{0}\|$, $\left\|\frac{\zeta_0}{\sqrt{\gamma-1}}\right\|$ and $\|\xi_{0}\|_1$, it holds for $(t,x)\in[0,T]\times\mathbb{R}$ that
\begin{eqnarray}\label{2.27}
C_0^{-1}\leq v(t,x)\leq C_0,
\end{eqnarray}
\end{Lemma}
\noindent{\bf Proof.}~~By the Cauchy inequality and $(\ref{2.8})$, we have
\begin{equation}\label{2.28}
C_8^{-1}\leq B(t,x)\leq C_8,~~~x\in[i,i+1],~t\geq 0,
\end{equation}
where $C_8^{-1},C_8$ are two positive constants only depending on $C_4$.

Integrating $(\ref{2.21})_4$ over $(s,t)\times [i+1,i+2]$, then applying the Cauchy inequality, the H\"{o}lder inequality and the Jessen's inequality for the function $\frac{1}{x}(x>0)$, using $(\ref{2.8})$ and Lemma 2.2, we obtain that
\begin{eqnarray}\label{2.29}
\int_s^t\int_{i+1}^{i+2}\mathcal{E}(\tau,x)dxd\tau=&&\int_s^t\int_{i+1}^{i+2}\left(\nu\frac{\psi_x}{v}+\nu\frac{U_x}{v}-\frac{\chi_x^2}{2v^2}-\frac{R\theta}{v}\right)dxd\tau\nonumber\\
\leq &&\int_s^t\int_{i+1}^{i+2}\left(\nu\frac{\psi_x}{v}+\nu\frac{U_x}{v}-\frac{R\theta}{v}\right)dxd\tau\nonumber\\
\leq &&C\int_s^t\int_{i+1}^{i+2}\frac{\nu\Theta}{v\theta}\psi_x^2dxd\tau-\frac{R}{2}\int_s^t\int_{i+1}^{i+2}\frac{\theta}{v}dxd\tau+\int_s^t\int_{i+1}^{i+2}\nu\frac{U_x}{v}dxd\tau\nonumber\\
\leq &&C_{9}-\frac{R\underline{\Theta}}{2}\int_s^t\int_{i+1}^{i+2}\frac{1}{v}dxd\tau+C(M_0^{-1})\left(\int_s^t\int_{i+1}^{i+2}U_x^2dxd\tau\right)^{\frac{1}{2}}\sqrt{t-s}\nonumber\\
\leq &&C_{9}-\frac{R\underline{\Theta}}{2\overline{V}}\int_s^t\left(\int_{i+1}^{i+2}vdx\right)^{-1}d\tau+C(M_0^{-1})\varepsilon^{\frac{1}{4}}(t-s)\\
\leq &&C_{9}-\frac{R\underline{\Theta}}{2\alpha_2\overline{V}}(t-s)+C(M_0^{-1})\varepsilon^{\frac{1}{4}}(t-s)\nonumber\\
\leq &&C_{9}-\frac{R\underline{\Theta}}{4\alpha_2\overline{V}}(t-s):=C_{9}-\frac{t-s}{C_{10}},\nonumber
\end{eqnarray}
where letting $C(M_0^{-1})\varepsilon^{\frac{1}{4}}=\frac{R\underline{\Theta}}{4\alpha_2\overline{V}}$, provided that $\varepsilon$ is sufficiently small, and $C_{10}=\left(\frac{R\underline{\Theta}}{4\alpha_2\overline{V}}\right)^{-1}$.

It follows from the definition of $Y(t)$ and $(\ref{2.29})$ that
\begin{equation}\label{2.30}
0\leq Y(t)\leq C_{11}e^{-t/C_{12}},~~\frac{Y(t)}{Y(s)}\leq C_{11}e^{-(t-s)/C_{12}},
\end{equation}
where $C_{11}:=e^{\frac{1}{\nu}C_9}$ and $C_{12}:=\nu C_{10}$.

Integrating (\ref{2.20}) over $[i+1,i+2]$ with respect to $x$, using Lemmas 4.1-4.2 and (\ref{2.30}), we obtain
\begin{eqnarray}\label{2.31}
\alpha_1\underline{V}\leq\int_{i+1}^{i+2}\widetilde{v}(t,x)dx\leq&& \overline{V}C_8Y(t)+\overline{V}C_8^2\int_0^t\frac{Y(t)}{Y(\tau)}\int_{i+1}^{i+2}\left(R\theta+\frac{\chi_x^2}{2v}\right)dxd\tau\\
\leq&&C_{11}\overline{V}C_8e^{-t/C_{12}}+2\overline{V}C_4C_8^2\int_0^t\frac{Y(t)}{Y(\tau)}d\tau.\nonumber
\end{eqnarray}
This directly yields that
\begin{equation}\label{2.32}
\int_0^t\frac{Y(t)}{Y(\tau)}d\tau\geq \frac{\alpha_1\underline{V}C_8^{-1}-C_{11}\overline{V}e^{-t/C_{12}}}{2\overline{V}C_8C_4}:=C_{13}-C_{14}e^{-t/C_{12}}.
\end{equation}
From $(\ref{2.6-1})_2$, (\ref{2.20}), (\ref{2.28}) and (\ref{2.32}), we have
\begin{eqnarray}\label{2.33}
v(t,x)\geq&& \frac{1}{\nu C_8^2}\int_0^t\frac{Y(t)}{Y(\tau)}\left(\theta+\frac{\chi_x^2}{v}\right)d\tau
\geq \frac{1}{\nu C_8^2}\int_0^t\frac{Y(t)}{Y(\tau)}\theta(\tau,x) d\tau\\
\geq&&\frac{1}{\nu C_8^2}\left(C_{13}\underline{\Theta}-C_{14}\underline{\Theta}e^{-t/C_{12}}\right)
\geq\frac{1}{\nu C_8^2}\frac{C_{13}\underline{\Theta}}{2}:=C_{15}, ~~~~~~~x\in\mathbb{R},~t\geq T_0,\nonumber
\end{eqnarray}
where $C_{13}$, $C_{14}$ are some positive constant depending on $C_4$, and $T_0$ and $C_{15}$ are positive constants independent of $t$.

Integrating (\ref{2.20}) over $[i+1,i+2]$ with respect to $x$, using Lemmas 4.1-4.2 and (\ref{2.30}), we get
\begin{eqnarray*}
\alpha_1\underline{V}Y(t)^{-1}\leq&& \overline{V}C_8+2C_4\overline{V}C_8^2\int_0^t\frac{1}{Y(\tau)}d\tau.
\end{eqnarray*}
Then the above inequality by applying the Gronwall inequality yields that
\begin{eqnarray*}
Y(t)^{-1}\leq&&\frac{\overline{V}C_8}{\alpha_1\underline{V}}\exp\left\{\frac{2C_4\overline{V}C_8^2}{\alpha_1\underline{V}}T_0\right\},~ t\in[0,T_0].
\end{eqnarray*}
That is,
\begin{eqnarray*}
Y(t)\geq&&\frac{\alpha_1\underline{V}}{\overline{V}C_8}\exp\left\{-2C_4\overline{V}\frac{C_8^2}{\alpha_1\underline{V}}T_0\right\}:=C_{17}e^{-C_{16}T_0},~ t\in[0,T_0]. \\
\end{eqnarray*}
Due to
\begin{equation}\label{2.34}
v(t,x)\geq B(t,x)Y(t)\geq C_8^{-1}C_{17}e^{-C_{16}T_0}=C_{18},~\forall t\in[0,T_0],~x\in\mathbb{R}.
\end{equation}
Combining with (\ref{2.33})-(\ref{2.34}), we obtain
\begin{equation}\label{2.35}
v(t,x)\geq \min\left\{C_{15},C_{18}\right\}:=C_0^{-1},~\forall t\in[0,+\infty),~x\in\mathbb{R}.
\end{equation}
Now, together (\ref{2.30}) with (\ref{2.20}) gives
\begin{eqnarray}\label{2.36}
v(t,x)\leq&& C_8C_{11}e^{-t/C_{12}}+C_{11}C_8^2\int_0^t\left(R\theta+\frac{\chi_x^2}{2v}\right)e^{-\frac{(t-\tau)}{C_{12}}}d\tau\\
\leq&& C_8C_{11}+C_{11}C_8^2\int_0^t\left(R\overline{\Theta}e^{-\frac{(t-\tau)}{C_{12}}}+\max_{x\in\mathbb{R}}\left(\frac{\chi_x^2}{2v^2}\right)\max_{x\in\mathbb{R}} v(t,x)e^{-\frac{(t-\tau)}{C_{12}}}\right)d\tau.\nonumber
\end{eqnarray}

On the other hand, for $\forall(t,x)\in[0,T]\times[i+1,i+2]$, there exists a positive constant $c_i(t)\in[i+1,i+2]$, such that
\begin{eqnarray}\label{2.37}
\frac{\chi_x^2}{v^2}(t,x)-\frac{\chi_x^2}{v^2}(t,c_i(t))=&&\int_{c_i(t)}^x2\frac{\chi_y}{v}\left(\frac{\chi_y}{v}\right)_ydy\nonumber\\
\leq&&C\int_{i+1}^{i+2}\left|\frac{\chi_x}{v}\right||-\mu+\chi(\chi-1)(\chi+1)|dx\\
\leq&&\epsilon\max_{x\in\mathbb{R}}\left(\frac{\chi_x}{v}\right)^2\int_{\mathbb{R}}\frac{\theta}{v}dx+C_\epsilon\int_{i+1}^{i+2}\frac{v}{\theta}\mu^2dx+C_\epsilon\nonumber\\
\leq&&C\epsilon\max_{x\in\mathbb{R}}\left(\frac{\chi_x}{v}\right)^2+C_\epsilon\int_{\mathbb{R}}\frac{v}{\theta}\mu^2dx+C_\epsilon,\nonumber
\end{eqnarray}
where we used the Cauchy inequality, taking $\epsilon$ suitably small, (\ref{2.35}), $(\ref{2.1})_{3,4}$, $(\ref{2.6-1})_2$ and (\ref{2.8}).\\
The above inequality implied that
\begin{eqnarray}\label{2.38}
\max_{x\in\mathbb{R}}\left(\frac{\chi_x^2}{v^2}\right)(t,x)\leq C\int_{\mathbb{R}}\frac{v}{\theta}\mu^2dx+C.
\end{eqnarray}

Substituting (\ref{2.38}) into (\ref{2.36}), and applying the Gronwall inequality and Lemma 4.1, we have
\begin{eqnarray}\label{2.39}
v(t,x)\leq C_0,~x\in\mathbb{R},~t\geq 0,
\end{eqnarray}
which, together with (\ref{2.35}), completes the proof of Lemma 4.4.

From Lemma 4.1, Lemma 4.4 and $(\ref{2.6-1})_2$, we have the following corollary.
\begin{Corollary} Under the assumptions of Proposition 3.2, then there exists a uniform positive constant $C_{19}>0$ depending only on $\underline{V},\overline{V}$, $\underline{\Theta}$, $\overline{\Theta}$, $\delta$, $\|\phi_{0}\|$, $\|\psi_{0}\|$, $\left\|\frac{\zeta_0}{\sqrt{\gamma-1}}\right\|$ and $\|\xi_{0}\|_1$, it holds for $(t,x)\in[0,T]\times\mathbb{R}$ that
\begin{eqnarray}\label{2.40}
     &&\left\|\left(\phi,\psi,\frac{\zeta}{\sqrt{\gamma-1}},\xi,\xi_x\right)(t)\right\|^2+\int_0^t\|(\mu,\psi_x,\zeta_x)(\tau)\|^2d\tau\nonumber\\
     \leq&&C_{19}\left(\left\|\left(\phi_0,\psi_0,\frac{\zeta_0}{\sqrt{\gamma-1}}\right)\right\|^2+\|\xi_{0}\|_1^2+1\right).
\end{eqnarray}
\end{Corollary}

\begin{Lemma} Under the assumptions of Proposition 3.2, then there exists a uniform positive constant $C_{20}>0$ depending only on $\underline{V},\overline{V}$, $\underline{\Theta}$, $\overline{\Theta}$, $\delta$, $\|\phi_{0}\|$, $\|\psi_{0}\|$, $\left\|\frac{\zeta_0}{\sqrt{\gamma-1}}\right\|$ and $\|\xi_{0}\|_1$, it holds for $(t,x)\in[0,T]\times\mathbb{R}$ that
\begin{eqnarray}\label{2.41}
\|\xi(t)\|^2+\int_0^t\|\xi(\tau)\|^2_1d\tau\leq C_{20}\left(\left\|\left(\phi_0,\psi_0,\frac{\zeta_0}{\sqrt{\gamma-1}}\right)\right\|^2+\|\xi_{0}\|_1^2+1\right).
\end{eqnarray}
\end{Lemma}
\noindent{\bf Proof.}~~Multiplying $(\ref{2.1})_{3.4}$ by $\xi$, and then integrating the resultant equation with respect to $t$ and $x$ over $[0,t]\times\mathbb{R}$, we obtain
\begin{eqnarray}\label{2.42}
&&\int_\mathbb{R}\frac{\xi^2}{2v}dx+\int_0^t\int_\mathbb{R}\left(\xi^2(\xi+1)(\xi+2)+\frac{\xi^2U_x}{2v^2}+\frac{\xi_x^2}{v}\right)dxd\tau\nonumber\\
=&&\int_\mathbb{R}\frac{\xi_0^2}{2v_0}dx-\int_0^t\int_\mathbb{R}\frac{\xi^2\psi_x}{2v^2}dxd\tau\nonumber\\
\leq&& C\|\xi_{0}\|^2+C\int_0^t\|\xi(\tau)\|^{\frac{3}{2}}\|\xi_x(\tau)\|^{\frac{1}{2}}\|\psi_x(\tau)\|d\tau\\
\leq&& C\|\xi_{0}\|^2+\eta\int_0^t\|\xi_x(\tau)\|^2d\tau+C_\eta\int_0^t\|\xi(\tau)\|^{2}\|\psi_x(\tau)\|^{\frac{4}{3}}d\tau\nonumber\\
\leq&& C\|\xi_{0}\|^2+\eta\int_0^t\|\xi_x(\tau)\|^2d\tau+C_\eta\left(\int_0^t\|\xi(\tau)\|^2d\tau\right)^{\frac{1}{3}}\left(\int_0^t\|\psi_x(\tau)\|^2d\tau\right)^{\frac{2}{3}}\nonumber\\
\leq&& C\|\xi_{0}\|^2+\eta\int_0^t\|\xi_x(\tau)\|^2d\tau+\frac{1}{3}\int_0^t\|\xi(\tau)\|^2d\tau+C_\eta\int_0^t\|\psi_x(\tau)\|^2d\tau,\nonumber
\end{eqnarray}
where we used to the H\"{o}lder inequality, the Cauchy inequality, the Young inequality and Lemma 4.4.

By applying Corollary 4.1 and taking $\eta$ suitably small, we get (\ref{2.41}) immediately. This complete the proof of Lemma 4.5.

\begin{Lemma} Under the assumptions of Proposition 3.2, then there exists a uniform positive constant $C_{21}>0$ depending only on $\underline{V},\overline{V}$, $\underline{\Theta}$, $\overline{\Theta}$, $\delta$, $\|\phi_{0}\|$, $\|\psi_{0}\|$, $\left\|\frac{\zeta_0}{\sqrt{\gamma-1}}\right\|$ and $\|\xi_{0}\|_1$, it holds for $(t,x)\in[0,T]\times\mathbb{R}$ that
\begin{eqnarray}\label{2.43}
&&\|\phi_x(t)\|^2+\int_0^t\|\phi_x(\tau)\|^2d\tau\\
\leq&& C_{21}\left(\left\|\left(\psi_0,\frac{\zeta_0}{\sqrt{\gamma-1}}\right)\right\|^2+\|(\phi_0,\xi_{0})\|_1^2+1+\int_0^t\|\xi_t(\tau)\|^2\|\mu(\tau)\|^2d\tau\right).\nonumber
\end{eqnarray}
\end{Lemma}
\noindent{\bf Proof.}~~Inserting $(\ref{2.1})_{1}$ into $(\ref{2.1})_{2}$, then we can rewritten $(\ref{2.1})_{2}$
\begin{eqnarray}\label{2.44}
\left(\nu\frac{\phi_x}{v}-\psi\right)_t+\frac{R\theta\phi_x}{v^2}=&&\frac{R\zeta_x}{v}-\frac{R\Theta_x\phi}{vV}+R\left(\frac{\Theta}{V^2}-\frac{\theta}{v^2}\right)V_x-\nu\frac{U_{xx}}{v}\nonumber\\
&&+\nu\frac{\psi_xV_x+U_xV_x}{v^2}+g(V,\Theta)_x+\frac{1}{2}\left(\frac{\chi_x^2}{v^2}\right)_x.
\end{eqnarray}
Multiplying (\ref{2.44}) by $\frac{\phi_x}{v}$, then we integrate the resultant equation with respect to $t$ and $x$  over $[0,t]\times\mathbb{R}$
\begin{eqnarray}\label{2.45}
&&\int_\mathbb{R}\nu\frac{\phi_x^2}{v^2}dx+\int_0^t\int_\mathbb{R}\frac{R\theta\phi_x^2}{v^3}dxd\tau\nonumber\\
=&&\int_\mathbb{R}\left(\nu\frac{\phi_{0x}^2}{v_0^2}dx-\psi_0\frac{\phi_{0x}}{v_0}\right)dx+\int_\mathbb{R}\psi\frac{\phi_{x}}{v}dx+\int_0^t\int_\mathbb{R}\frac{\psi_{x}^2}{v}dxd\tau\nonumber\\
&&+\int_0^t\int_\mathbb{R}\frac{\phi_{x}}{v}\left\{\frac{R\zeta_x}{v}-\frac{R\Theta_x\phi}{vV}+R\left(\frac{\Theta}{V^2}-\frac{\theta}{v^2}\right)V_x\right\}dxd\tau\\
&&+\int_0^t\int_\mathbb{R}\left(\frac{\psi}{v^2}(\phi_xU_x-\psi_xV_x)+\nu\frac{(\psi_xV_x+U_xV_x)\phi_x}{v^3}-\nu\frac{U_{xx}\phi_x}{v^2}\right)dxd\tau\nonumber\\
&&+\int_0^t\int_\mathbb{R}g(V,\Theta)_x\frac{\phi_{x}}{v}dxd\tau+\frac{1}{2}\int_0^t\int_\mathbb{R}\left(\frac{\chi_x^2}{v^2}\right)_x\frac{\phi_{x}}{v}dxd\tau\nonumber\\
=&&\int_\mathbb{R}\left(\nu\frac{\phi_{0x}^2}{v_0^2}dx-\psi_0\frac{\phi_{0x}}{v_0}\right)dx+\int_\mathbb{R}\psi\frac{\phi_{x}}{v}dx+\int_0^t\int_\mathbb{R}\frac{\psi_{x}^2}{v}dxd\tau+\sum_{i=1}^4\mathcal{I}_i.\nonumber
\end{eqnarray}

Now we deal with the second and third terms on the right-hand side of (\ref{2.45}), and $\mathcal{I}_i(i=1,2,3,4)$ term by term. For this purpose, applying the Cauchy inequality, the Young inequality, $(\ref{2.6-1})_2$, Lemma 2.2, Lemmas 4.4-4.5 and Corollary 4.1, we get
\begin{eqnarray}\label{2.46}
\int_\mathbb{R}\psi\frac{\phi_{x}}{v}dx+\int_0^t\int_\mathbb{R}\frac{\psi_{x}^2}{v}dxd\tau\leq \frac{1}{4}\int_\mathbb{R}\nu\frac{\phi_x^2}{v^2}dx+C\|\psi(t)\|^2+C\int_0^t\|\psi_x(\tau)\|^2d\tau,
\end{eqnarray}
\begin{eqnarray}\label{2.47}
|\mathcal{I}_1|
\leq &&\frac{1}{8}\int_0^t\int_\mathbb{R}\frac{R\theta\phi_x^2}{v^3}dxd\tau+C\int_0^t\left(\|\zeta_x(\tau)\|^2+\left\|\left(\frac{\zeta}{\sqrt{\gamma-1}},\phi\right)(\tau)\right\|^2\|(V_x,\Theta_x)(\tau)\|^2_{L^\infty}\right)d\tau\nonumber\\
\leq &&\frac{1}{8}\int_0^t\int_\mathbb{R}\frac{R\theta\phi_x^2}{v^3}dxd\tau+C\int_0^t\|\zeta_x(\tau)\|^2d\tau+C,
\end{eqnarray}
\begin{eqnarray}\label{2.48}
&&|\mathcal{I}_2|+|\mathcal{I}_3|\nonumber\\
\leq &&\frac{1}{8}\int_0^t\int_\mathbb{R}\frac{R\theta\phi_x^2}{v^3}dxd\tau+C\int_0^t\int_\mathbb{R}\left(\psi_x^2V_x^2+\psi\psi_xV_x+\psi^2U_x^2+U_x^2V_x^2+U_{xx}^2+g(V,\Theta)_x^2\right)dxd\tau\nonumber\\
\leq &&\frac{1}{8}\int_0^t\int_\mathbb{R}\frac{R\theta\phi_x^2}{v^3}dxd\tau+C\int_0^t\left(\|\psi(\tau)\|^2\|(V_x,U_x)(\tau)\|^2_{L^\infty}+\|(\psi_x,U_{xx})(\tau)\|^2+\|g(V,\Theta)_x(\tau)\|^2\right)d\tau\nonumber\\
\leq &&\frac{1}{8}\int_0^t\int_\mathbb{R}\frac{R\theta\phi_x^2}{v^3}dxd\tau+C\int_0^t\|\psi_x(\tau)\|^2d\tau+C,
\end{eqnarray}
and
\begin{eqnarray}\label{2.49}
|\mathcal{I}_4|\leq &&\frac{1}{8}\int_0^t\int_\mathbb{R}\frac{R\theta\phi_x^2}{v^3}dxd\tau+C\int_0^t\int_\mathbb{R}\left(\frac{\xi_x}{v}\right)^2\left[\left(\frac{\xi_x}{v}\right)_x\right]^2dxd\tau\nonumber\\
\leq &&\frac{1}{8}\int_0^t\int_\mathbb{R}\frac{R\theta\phi_x^2}{v^3}dxd\tau+C\int_0^t\|\xi_x(\tau)\|\|\xi_{xx}(\tau)\|\left(\|\mu(\tau)\|^2+\|(\chi^2-1)(\tau)\|^2\right)d\tau\\
\leq &&\frac{1}{8}\int_0^t\int_\mathbb{R}\frac{R\theta\phi_x^2}{v^3}dxd\tau+C\int_0^t\|\xi_x(\tau)\|\left(\|\mu(\tau)\|+\|\xi_x(\tau)\|(\|\phi_x(\tau)\|^2+\|V_x(\tau)\|^2)+\|(\chi^2-1)(\tau)\|\right)\nonumber\\
&&\cdot\left(\|\mu(\tau)\|^2+\|(\chi^2-1)(\tau)\|^2\right)d\tau\nonumber\\
\leq &&\frac{1}{8}\int_0^t\int_\mathbb{R}\frac{R\theta\phi_x^2}{v^3}dxd\tau+C\int_0^t\left(\|\phi_x(\tau)\|^2\|(\xi_t,\xi_x)(\tau)\|^2+\|\xi_t(\tau)\|^2\|\mu(\tau)\|^2+\|\mu(\tau)\|^2+\|\xi(\tau)\|^2_1\right)d\tau,\nonumber
\end{eqnarray}
where we used the following inequality:
\begin{eqnarray*}
\|\xi_{xx}(t)\|\leq&& C(\|\mu(t)\|+\|\xi_x(\phi_x+V_x)(t)\|+\|(\chi^2-1)(t)\|)\nonumber\\
\leq&& C(\|\mu(t)\|+\|\xi_x(t)\|^\frac{1}{2}\|\xi_{xx}(t)\|^\frac{1}{2}(\|\phi_x(t)\|+\|V_x(t)\|)+\|(\chi^2-1)(t)\|)\\
\leq&&\frac{1}{2}\|\xi_{xx}(t)\|+C(\|\mu(t)\|+\|\xi_x(t)\|(\|\phi_x(t)\|^2+\|V_x(t)\|^2)+\|(\chi^2-1)(t)\|),\nonumber
\end{eqnarray*}
that is, \begin{eqnarray}\label{2.50}
\|\xi_{xx}(t)\|\leq C\biggl(\|\mu(t)\|+\|\xi_x(t)\|(\|\phi_x(t)\|^2+\|V_x(t)\|^2)+\|(\chi^2-1)(t)\|\biggr).
\end{eqnarray}
Then inserting the estimates (\ref{2.46})-(\ref{2.49}) into the identity for (\ref{2.45}) and using the Gronwall inequality, Lemmas 4.4-4.5 and Corollary 4.1, it yields that (\ref{2.43}). This complete the proof of Lemma 4.6.

\begin{Lemma} Under the assumptions of Proposition 3.2, then there exists a uniform positive constant $C_{22}>0$ depending only on $\underline{V},\overline{V}$, $\underline{\Theta}$, $\overline{\Theta}$, $\delta$, $\|\phi_{0}\|$, $\|\psi_{0}\|$, $\left\|\frac{\zeta_0}{\sqrt{\gamma-1}}\right\|$ and $\|\xi_{0}\|_1$, it holds for $(t,x)\in[0,T]\times\mathbb{R}$ that
\begin{eqnarray}\label{2.51}
&&\|\xi_t(t)\|^2+\int_0^t\|\xi_t(\tau)\|^2_1d\tau\\
\leq && C_{22}\left(\left\|\left(\psi_0,\frac{\zeta_0}{\sqrt{\gamma-1}}\right)\right\|^2+\|\phi_0\|^2_1+\|\xi_{0}\|_2^2+1+\int_0^t\|\psi_x(\tau)\|^2\|\phi_x(\tau)\|^2d\tau\right).\nonumber
\end{eqnarray}
\end{Lemma}
\noindent{\bf Proof.}~~Let us rewritten $(\ref{2.1})_{3,4}$ as follows:
\begin{equation}\label{2.52}
\frac{\xi_t}{v}=\left(\frac{\xi_x}{v}\right)_x-\xi(\xi+1)(\xi+2).
\end{equation}
Taking $(\ref{2.52})_t\times\xi_t$, and integrating the resultant equation with respect to $t$ and $x$ over $[0,t]\times\mathbb{R}$, we obtain
\begin{eqnarray}\label{2.53}
&&\frac{1}{2}\int_\mathbb{R}\frac{\xi_t^2}{v}dx+\int_0^t\int_\mathbb{R}\left((3\xi^2+2)\xi_t^2+\frac{\xi_{xt}^2}{v}\right)dxd\tau\nonumber\\
=&&\frac{1}{2}\int_\mathbb{R}\frac{\xi_{0t}^2}{v_0}dx+\int_0^t\int_\mathbb{R}\left(\frac{\xi_{t}^2(\psi_x+U_x)}{2v^2}+\frac{\xi_{x}\xi_{xt}(\psi_x+U_x)}{v^2}-6\xi\xi_t^2\right)dxd\tau\nonumber\\
\leq&&\frac{1}{2}\int_\mathbb{R}\frac{\xi_{0t}^2}{v_0}dx+\frac{1}{4}\int_0^t\int_\mathbb{R}\frac{\xi_{xt}^2}{v}dxd\tau+C\int_0^t\|(\xi_{t},\xi_x)(\tau)\|^2d\tau\nonumber\\
&&+C\int_0^t\left(\|\xi_t(\tau)\|^\frac{3}{2}\|\xi_{xt}(\tau)\|^\frac{1}{2}\|\psi_{x}(\tau)\|+\|\xi_{x}(\tau)\|\|\xi_{xx}(\tau)\|\|\psi_{x}(\tau)\|^2\right)d\tau\nonumber\\
\leq&&\frac{1}{2}\int_\mathbb{R}\frac{\xi_{0t}^2}{v_0}dx+\frac{1}{2}\int_0^t\int_\mathbb{R}\frac{\xi_{xt}^2}{v}dxd\tau+C\int_0^t\left[\|\xi_{t}(\tau)\|^2(\|\psi_{x}(\tau)\|^2+1)+\|\xi_x(\tau)\|^2\right]d\tau\\
&&+C\int_0^t\|\xi_{x}(\tau)\|\biggl(\|\mu(\tau)\|+\|\xi_x(\tau)\|(\|\phi_x(\tau)\|^2+\|V_x(\tau)\|^2)+\|(\chi^2-1)(\tau)\|\biggr)\|\psi_{x}(\tau)\|^2d\tau\nonumber\\
\leq&&\frac{1}{2}\int_\mathbb{R}\frac{\xi_{0t}^2}{v_0}dx+\frac{1}{2}\int_0^t\int_\mathbb{R}\frac{\xi_{xt}^2}{v}dxd\tau+C\int_0^t\|\xi_{t}(\tau)\|^2\|\psi_{x}(\tau)\|^2d\tau\nonumber\\
&&+C\int_0^t\|(\xi_t,\psi_x,\xi_x)(\tau)\|^2d\tau+C\int_0^t\|\psi_x(\tau)\|^2\|\phi_x(\tau)\|^2d\tau,\nonumber
\end{eqnarray}
where used the Cauchy inequality, the Young inequality, Lemma 2.2, Lemmas 4.4-4.5, Corollary 4.1 and (\ref{2.50}).

Then (\ref{2.51}) follows immediately from (\ref{2.53}), the Gronwall inequality, Lemmas 4.4-4.5, Corollary 4.1 and the fact that $$\|\xi_{0t}\|^2=\left\|\displaystyle \xi_{0xx}-\frac{\xi_{0x}v_{0x}}{v_0}-v_0\xi_0(\xi_0+1)(\xi_0+2)\right\|^2\leq C\left(\|(\xi_{0xx},\xi_{0x},\xi_{0}, \phi_{0x})\|^2+1\right).$$

This complete the proof of Lemma 4.7.

\begin{Corollary} Under the assumptions of Proposition 3.2, then there exist two uniform positive constant $C_{23}>0$ and $C_{24}>0$ depending only on $\underline{V},\overline{V}$, $\underline{\Theta}$, $\overline{\Theta}$, $\delta$, $\|\phi_{0}\|$, $\|\psi_{0}\|$, $\left\|\frac{\zeta_0}{\sqrt{\gamma-1}}\right\|$ and $\|\xi_{0}\|_1$, it holds for $(t,x)\in[0,T]\times\mathbb{R}$ that
\begin{eqnarray}\label{2.54}
\left\|\left(\phi_x,\xi_t\right)(t)\right\|^2+\int_0^t\|(\xi_t,\xi_{tx},\phi_x)(\tau)\|^2d\tau
     \leq C_{23}\left(\left\|\left(\psi_0,\frac{\zeta_0}{\sqrt{\gamma-1}}\right)\right\|^2+\|\phi_0\|^2_1+\|\xi_{0}\|_2^2+1\right),
\end{eqnarray}
and
\begin{eqnarray}\label{2.55}
    \left\|\xi_{xx}(t)\right\|^2+\int_0^t\|\xi_{xx}(\tau)\|^2d\tau
     \leq C_{24}\left(\left\|\left(\psi_0,\frac{\zeta_0}{\sqrt{\gamma-1}}\right)\right\|^2+\|\phi_0\|^2_1+\|\xi_{0}\|_2^2+1\right).
\end{eqnarray}
\end{Corollary}
\noindent{\bf Proof.}~~From Lemmas 4.6-4.7, we know
\begin{eqnarray*}
&&\left\|\left(\phi_x,\xi_t\right)(t)\right\|^2+\int_0^t\|(\xi_t,\xi_{tx},\phi_x)(\tau)\|^2d\tau\\
     \leq && \left(\left\|\left(\psi_0,\frac{\zeta_0}{\sqrt{\gamma-1}}\right)\right\|^2+\|\phi_0\|^2_1+\|\xi_{0}\|_2^2+1+\int_0^t\|\xi_t(\tau)\|^2\|\mu(\tau)\|^2d\tau+\int_0^t\|\psi_x(\tau)\|^2\|\phi_x(\tau)\|^2d\tau\right),\nonumber
\end{eqnarray*}
Then applying the Gronwall inequality and Corollary 4.1,  we get (\ref{2.54}) immediately.
After integrating (\ref{2.50}) over $[0,t]$, and then adding (\ref{2.50}), from (\ref{2.54}), Corollary 4.1 and Lemma 4.5,  we can obtain (\ref{2.55}) immediately.

\begin{Lemma} Under the assumptions of Proposition 3.2, then there exists a uniform positive constant $C_{25}>0$ depending only on $\underline{V},\overline{V}$, $\underline{\Theta}$, $\overline{\Theta}$, $\delta$, $\|\phi_{0}\|_1$, $\|\psi_{0}\|$, $\left\|\frac{\zeta_0}{\sqrt{\gamma-1}}\right\|$ and $\|\xi_{0}\|_2$, it holds for $(t,x)\in[0,T]\times\mathbb{R}$ that
\begin{eqnarray}\label{2.56}
\|\psi_x(t)\|^2+\int_0^t\|\psi_{xx}(\tau)\|^2d\tau
\leq C_{25}\left(\|(\phi_0,\psi_0)\|^2_1+\left\|\frac{\zeta_0}{\sqrt{\gamma-1}}\right\|^2+\|\xi_{0}\|_2^2+1\right).
\end{eqnarray}
\end{Lemma}
\noindent{\bf Proof.}~~Differentiating $(\ref{2.1})_2$ with respect to $x$ once, and then multiplying the resultant equation by $\psi_{x}$, we have
\begin{eqnarray}\label{2.57}
\left(\frac{\psi_x^2}{2}\right)_t+\nu\frac{\psi_{xx}^2}{v}=&&(\psi_t\psi_x)_x+\left(\frac{R\zeta_x}{v}-\frac{R\theta\phi_x}{v^2}-\frac{p_V(V,\Theta)V_x+p_\Theta(V,\Theta)\Theta_x}{v}\phi-\frac{R\zeta-p(V,\Theta)\phi}{v^2}V_x\right)\psi_{xx}\nonumber\\
&&+\left(-\nu\frac{U_{xx}}{v}+\nu\frac{(\psi_x+U_x)(\phi_x+V_x)}{v^2}+g(V,\Theta)_x+\frac{1}{2}\left(\frac{\xi_x^2}{v^2}\right)_x\right)\psi_{xx}\\
:=&&(\psi_t\psi_x)_x+\mathcal{R}_5+\mathcal{R}_6.\nonumber
\end{eqnarray}

We integrate (\ref{2.57}) with respect to $t$ and $x$ over $[0,t]\times\mathbb{R}$ to obtain
\begin{eqnarray}\label{2.58}
\frac{1}{2}\int_\mathbb{R}\psi_x^2dx+\int_0^t\int_\mathbb{R}\nu\frac{\psi_{xx}^2}{v}dxd\tau=\frac{1}{2}\int_\mathbb{R}\psi_{0x}^2dx+\int_0^t\int_\mathbb{R}\mathcal{R}_5dxd\tau+\int_0^t\int_\mathbb{R}\mathcal{R}_6dxd\tau.
\end{eqnarray}

Using the Cauchy inequality, the Young inequality, Lemma 2.2, Lemmas 4.4-4.5, Corollaries 4.1-4.2, one has that
\begin{eqnarray}\label{2.59}
&&\int_0^t\int_\mathbb{R}\mathcal{R}_5dxd\tau\nonumber\\
\leq && \frac{1}{4}\int_0^t\int_\mathbb{R}\nu\frac{\psi_{xx}^2}{v}dxd\tau+C\int_0^t\int_\mathbb{R}\left[\zeta_x^2+\phi_x^2+\phi^2(V_x^2+\Theta_x^2)+(\zeta^2+\phi^2)V_x^2\right]dxd\tau\\
\leq &&\frac{1}{4}\int_0^t\int_\mathbb{R}\nu\frac{\psi_{xx}^2}{v}dxd\tau+C\int_0^t\left(\|(\phi_x,\zeta_x)(\tau)\|^2+\left\|\left(\phi,\frac{\zeta}{\sqrt{\gamma-1}}\right)(\tau)\right\|^2\|(V_x,\Theta_x)(\tau)\|^2_{L^\infty}\right)d\tau\nonumber\\
\leq &&\frac{1}{4}\int_0^t\int_\mathbb{R}\nu\frac{\psi_{xx}^2}{v}dxd\tau+C\int_0^t\|(\phi_x,\zeta_x)(\tau)\|^2d\tau+C,\nonumber
\end{eqnarray}
and
\begin{eqnarray}\label{2.60}
\int_0^t\int_\mathbb{R}\mathcal{R}_6dxd\tau\leq &&\frac{1}{4}\int_0^t\int_\mathbb{R}\nu\frac{\psi_{xx}^2}{v}dxd\tau+C\int_0^t\biggl(\|(U_{xx},g(V,\Theta)_x,\psi_x,\phi_x)(\tau)\|^2\nonumber\\
&&+\|\psi_x(\tau)\|\|\psi_{xx}(\tau)\|\|\phi_x(\tau)\|^2+\|\xi_x(\tau)\|^2_{L^\infty}\|(\xi_t,\xi)(\tau)\|^2\biggr)d\tau\\
\leq &&\frac{1}{2}\int_0^t\int_\mathbb{R}\nu\frac{\psi_{xx}^2}{v}dxd\tau+C\int_0^t\|(\xi_t,\xi_{x},\xi_{xx},\psi_x,\phi_x)(\tau)\|^2d\tau+C.\nonumber
\end{eqnarray}

Inserting (\ref{2.59}) and (\ref{2.60}) into (\ref{2.58}), applying Lemmas 4.4-4.5 and Corollaries 4.1-4.2, we get (\ref{2.56}) immediately. This complete the proof of the Lemma 4.8.

\begin{Lemma} Under the assumptions of Proposition 3.2, then there exists a uniform positive constant $C_{26}>0$ depending only on $\underline{V},\overline{V}$, $\underline{\Theta}$, $\overline{\Theta}$, $\delta$, $\|\phi_{0}\|_1$, $\|\psi_{0}\|_1$, $\left\|\frac{\zeta_0}{\sqrt{\gamma-1}}\right\|$ and $\|\xi_{0}\|_2$, it holds for $(t,x)\in[0,T]\times\mathbb{R}$ that
\begin{eqnarray}\label{2.61}
\left\|\frac{\zeta_x}{\sqrt{\gamma-1}}(t)\right\|^2+\int_0^t\|\zeta_{xx}(\tau)\|^2d\tau
\leq C_{26}\left(\left\|\left(\phi_0,\psi_0,\frac{\zeta_0}{\sqrt{\gamma-1}}\right)\right\|^2_1+\|\xi_{0}\|_2^2+1\right).
\end{eqnarray}
\end{Lemma}
\noindent{\bf Proof.}~~Differentiating $(\ref{2.1})_6$ with respect to $x$ once, and then multiplying the resultant equation by $\zeta_{x}$, we have
\begin{eqnarray}\label{2.62}
\left(\frac{C_v}{2}\zeta_x^2\right)_t+\kappa\frac{\zeta_{xx}^2}{v}=&&(C_v\zeta_t\zeta_x)_x+\left(-\frac{R\Theta_{xx}}{v}+\kappa\frac{(\zeta_x+\Theta_x)(\phi_x+V_x)}{v^2}\right)\zeta_{xx}\nonumber\\
&&-\left(\nu\frac{\psi_x^2+U_x^2+2\psi_xU_x}{v}+v\mu^2-r(V,\Theta)\right)\zeta_{xx}\\
:=&&(C_v\zeta_t\zeta_x)_x+\mathcal{R}_7+\mathcal{R}_8.\nonumber
\end{eqnarray}

Integrating (\ref{2.62}) with respect to $t$ and $x$ over $[0,t]\times\mathbb{R}$, then applying the Cauchy inequality, the Young inequality, Lemma 2.2, Lemma 4.4 and Corollaries 4.1-4.2, we get that
\begin{eqnarray}\label{2.63}
&&\frac{C_v}{2}\int_\mathbb{R}\zeta_x^2dx+\int_0^t\int_\mathbb{R}\kappa\frac{\zeta_{xx}^2}{v}dxd\tau\nonumber\\
=&&\frac{C_v}{2}\int_\mathbb{R}\zeta_{0x}^2dx+\int_0^t\int_\mathbb{R}\mathcal{R}_7dxd\tau+\int_0^t\int_\mathbb{R}\mathcal{R}_8dxd\tau\\
\leq&&C\|\zeta_{0x}\|^2+\frac{1}{4}\int_0^t\int_\mathbb{R}\kappa\frac{\zeta_{xx}^2}{v}dxd\tau+C\int_0^t\biggl(\|\Theta_{xx}(\tau)\|^2+\|\zeta_x(\tau)\|\|\zeta_{xx}(\tau)\|\|\phi_x(\tau)\|^2\nonumber\\
&&+\|(\zeta_x,\phi_x,\psi_x)(\tau)\|^2\|(V_x,\Theta_x,U_x)(\tau)\|^2_{L^\infty}+\|\psi_x(\tau)\|^3\|\psi_{xx}(\tau)\|+\|\xi_t(\tau)\|^3\|\xi_{tx}(\tau)\|\nonumber\\
&&+\|(\Theta_{x},U_x)(\tau)\|^2_{L^\infty}\|(V_x,U_x)(\tau)\|^2+\|r(V,\Theta)(\tau)\|^2\biggr)d\tau\nonumber\\
\leq&&C\|\zeta_{0x}\|^2+\frac{1}{2}\int_0^t\int_\mathbb{R}\kappa\frac{\zeta_{xx}^2}{v}dxd\tau+C\int_0^t\|(\phi_x,\psi_x,\psi_{xx},\zeta_x,\xi_t,\xi_{tx})(\tau)\|^2d\tau+C.\nonumber
\end{eqnarray}

One can obtain by using Lemmas 4.4-4.5, Lemma 4.8 and Corollaries 4.1-4.2 that get the estimate (\ref{2.61}) immediately. This complete the proof of the Lemma 4.9.

For the estimates of the second order derivatives of solutions $(\phi,\psi,\zeta,\varphi,\xi)(t,x)$ to the Cauchy problem (\ref{2.1})-(\ref{2.2}), by repeating the same argument as above, we can also obtain
\begin{Corollary} Under the assumptions of Proposition 3.2, then there exists a uniform positive constant $C_{27}>0$ depending only on $\underline{V},\overline{V}$, $\underline{\Theta}$, $\overline{\Theta}$, $\delta$, $\|\phi_{0}\|_2$, $\|\psi_{0}\|_2$, $\left\|\frac{\zeta_0}{\sqrt{\gamma-1}}\right\|_2$ and $\|\xi_{0}\|_2$, it holds for $(t,x)\in[0,T]\times\mathbb{R}$ that
\begin{eqnarray}\label{2.64}
&&\left\|\left(\phi_{xx},\psi_{xx},\zeta_{xx},\xi_{xx}\right)(t)\right\|^2+\int_0^t\left(\|\phi_{xx}(\tau)\|^2+\|(\psi_{xxx},\zeta_{xxx},\xi_{xxx})(\tau)\|^2\right)d\tau\nonumber\\
     \leq&& C_{27}\left(\left\|\left(\phi_0,\psi_0,\frac{\zeta_0}{\sqrt{\gamma-1}},\xi_{0}\right)\right\|_2^2+1\right).
\end{eqnarray}
\end{Corollary}

Now we are in a position to prove Proposition 3.2.\\
\noindent{\bf Proof of Proposition 3.2.} If $\varepsilon>0$ and $\gamma-1>0$ are sufficiently small such that $(\ref{2.6-1})_2$ and (\ref{2.16-1}) holds. Consequently, the conclusions of the above lemmas and corollaries hold.  Then we can obtain $(\ref{2.6})$-(\ref{2.7}) immediately.

\begin{center}
{\bf Acknowledgement}
\end{center}
This work is supported by the National
Natural Science Foundation of China under contracts 12171001, the Support
Program for Outstanding Young Talents in Universities of Anhui Province under contract gxyqZD2022007 and  the Excellent University Research and Innovation Team in Anhui Province under contract  2024AH010002.

\end{document}